\documentclass{amsart}
\def\Defined{}
\ifx\UseRussian\Defined
	\usepackage[T2A,T2B]{fontenc}
	\usepackage[cp1251]{inputenc}
	\usepackage[english,russian]{babel}
\fi
\usepackage{footmisc}
\usepackage[all]{xy}
\usepackage{color}
\definecolor{UrlColor}{rgb}{.9,0,.3}
\definecolor{SymbColor}{rgb}{.4,0,.9}
\definecolor{IndexColor}{rgb}{1,.3,.6}
\definecolor{eml1}{rgb}{.8,.1,.1}
\definecolor{eml2}{rgb}{.1,.6,.6}

\usepackage{xr-hyper}
\usepackage[unicode]{hyperref}
\hypersetup{pdfdisplaydoctitle=true}
\hypersetup{colorlinks}
\hypersetup{citecolor=UrlColor}
\hypersetup{urlcolor=UrlColor}
\hypersetup{pdffitwindow=true}
\hypersetup{pdfnewwindow=true}
\hypersetup{pdfstartview={FitH}}

\def\hyph{\penalty0\hskip0pt\relax-\penalty0\hskip0pt\relax}
\def\Hyph{-\penalty0\hskip0pt\relax}

\newcommand{\Basis}[1]{\overline{\overline{#1}}{}}
\newcommand{\Vector}[1]{\overline{#1}{}}
\newcommand{\gi}[1]{\boldsymbol{\textcolor{IndexColor}{#1}}}
\makeatletter
\newcommand{\NameDef}[1]{%
	\expandafter\gdef\csname #1\endcsname%
}%
\newcommand{\ShowSymbol}[1]{%
	\@nameuse{ViewSymbol#1}%
}%
\newcommand{\symb}[3]{%
	\@ifundefined{ViewSymbol#3}{%
		\NameDef{ViewSymbol#3}{\textcolor{SymbColor}{#1}}%
		\NameDef{RefSymbol#3}{\pageref{symbol: #3}}%
		\@namedef{LabelSymbol#3}{\label{symbol: #3}}%
	}{%
		\NameDef{RefSymbol#3}{}%
		\@namedef{LabelSymbol#3}{}%
	}%
	\ifcase#2
	\or
		$\@nameuse{ViewSymbol#3}$%
	\or
		\[\@nameuse{ViewSymbol#3}\]%
	\else%
	\fi%
	\@nameuse{LabelSymbol#3}%
}%
\makeatother

\newcommand{\subs}{${}_*$\Hyph}
\newcommand{\sups}{${}^*$\Hyph}

\newcommand{\CRstar}{{}_*{}^*}
\newcommand{\RCstar}{{}^*{}_*}

\newcommand{\RC}{$\RCstar$\Hyph}
\newcommand{\CR}{$\CRstar$\Hyph}
\newcommand{\drc}{$D\RCstar$\Hyph}
\newcommand{\dcr}{$D\CRstar$\hyph}
\newcommand{\rcd}{$\RCstar D$\Hyph}
\newcommand{\crd}{$\CRstar D$\Hyph}

\newcommand\sT{$\star T$\Hyph}%
\newcommand\Ts{$T\star$\Hyph}%

\renewcommand{\uppercasenonmath}[1]{}

\makeatletter
\newcommand\@dotsep{4.5}
\def\@tocline#1#2#3#4#5#6#7
{\relax
		\par \addpenalty\@secpenalty\addvspace{#2}%
		\begingroup \hyphenpenalty\@M
		\@ifempty{#4}{%
			\@tempdima\csname r@tocindent\number#1\endcsname\relax
		}{%
			\@tempdima#4\relax
		}%
		\parindent\z@ \leftskip#3\relax \advance\leftskip\@tempdima\relax
		\rightskip\@pnumwidth plus1em \parfillskip-\@pnumwidth
		#5\leavevmode\hskip-\@tempdima #6\relax
		\leaders\hbox{$\m@th
		\mkern \@dotsep mu\hbox{.}\mkern \@dotsep mu$}\hfill
		\hbox to\@pnumwidth{\@tocpagenum{#7}}\par
		\nobreak
		\endgroup
}
\makeatother 

\ifx\PrintBook\undefined
	\usepackage{fancyhdr}
	\makeatletter
	\renewcommand{\@indextitlestyle}{%
		\twocolumn[\section{\indexname}]%
		\def\IndexSpace{off}%
	}
	\makeatother 
	\thanks{\href{mailto:Aleks\_Kleyn@MailAPS.org}{Aleks\_Kleyn@MailAPS.org}}
\else
	\makeatletter
	\renewcommand{\@indextitlestyle}{%
		\twocolumn[\chapter{\indexname}]%
		\def\IndexSpace{off}%
		\let\@secnumber\@empty
		\chaptermark{\indexname}%
	}
	\makeatother 
	\email{\href{mailto:Aleks\_Kleyn@MailAPS.org}{Aleks\_Kleyn@MailAPS.org}}
\fi

\ifx\SelectlEnglish\undefined
	\ifx\UseRussian\undefined
		\def\SelectlEnglish{}
	\fi
\fi

\ifx\SelectlEnglish\undefined
	\newcommand\CurrentLanguage{Russian.}%
	\author{Александр Клейн}
	\newtheorem{theorem}{Теорема}[section]
	\newtheorem{corollary}[theorem]{Следствие}
	\theoremstyle{definition}
	\newtheorem{definition}[theorem]{Определение}
	\newtheorem{example}[theorem]{Пример}
	\newtheorem{xca}[theorem]{Exercise}
	\theoremstyle{remark}
	\newtheorem{remark}[theorem]{Замечание}
	
	\newcommand\Gbasis{$G$\Hyph базис}
	\newcommand\Gcoords{$G$\Hyph координат}
	\newcommand\Gspace{$G$\Hyph пространств}
	\newcommand\xRefDef[2]
		{\externaldocument[#1-Russian-]{#1.Russian}[http://arxiv.org/PS_cache/#2.pdf]}
	\newcommand\xRef[2]%
	{%
		\ifx\PrintBook\undefined%
		\citeBib{#1}-\ref{#1-Russian-#2}%
		\else%
		\ref{#2}%
		\fi%
	}%
	\newcommand\xEqRef[2]%
	{%
		\ifx\PrintBook\undefined%
		\citeBib{#1}-\eqref{#1-Russian-#2}%
		\else%
		\eqref{#2}%
		\fi%
	}%
	\ifx\PrintBook\undefined
		\newcommand{\BibTitle}{%
			\section{Список литературы}%
		}
	\else
		\newcommand{\BibTitle}{%
			\chapter{Список литературы}%
		}
	\fi
\else
	\newcommand\CurrentLanguage{English.}%
	\author{Aleks Kleyn}
	\newtheorem{theorem}{Theorem}[section]
	
	\theoremstyle{definition}
	\newtheorem{definition}[theorem]{Definition}
	\newtheorem{example}[theorem]{Example}
	
	\theoremstyle{remark}
	\newtheorem{remark}[theorem]{Remark}
	\newcommand\Gbasis{$G$\Hyph basis}
	\newcommand\Gcoords{$G$\Hyph coordinates}
	\newcommand\Gspace{$G$\Hyph space}
	\newcommand\xRefDef[2]
		{\externaldocument[#1-English-]{#1.English}[http://arxiv.org/PS_cache/#2.pdf]}
	\newcommand\xRef[2]%
	{%
		\ifx\PrintBook\undefined%
		\citeBib{#1}-\ref{#1-English-#2}%
		\else%
		\ref{#2}%
		\fi%
	}%
	\newcommand\xEqRef[2]%
	{%
		\ifx\PrintBook\undefined%
		\citeBib{#1}-\eqref{#1-English-#2}%
		\else%
		\eqref{#2}%
		\fi%
	}%
	\ifx\PrintBook\undefined
		\newcommand{\BibTitle}{%
			\section{References}%
		}
	\else
		\newcommand{\BibTitle}{%
			\chapter{References}%
		}
	\fi
\fi

\ifx\PrintBook\undefined
\else
	\numberwithin{section}{chapter}
\fi

\numberwithin{equation}{section}
\numberwithin{figure}{section}
\numberwithin{table}{section}
\numberwithin{Item}{section}
\numberwithin{Hfootnote}{section}

\makeatletter
\newcommand\org@maketitle{}
\let\org@maketitle\maketitle
\def\maketitle{%
	\hypersetup{pdftitle={\@title}}%
	\hypersetup{pdfauthor={\authors}}%
	\hypersetup{pdfsubject=\@keywords}%
	\org@maketitle
}
\def\make@stripped@name#1{%
	\begingroup
		\escapechar\m@ne
		\global\let\newname\@empty
		\protected@edef\Hy@tempa{\CurrentLanguage #1}%
		\edef\@tempb{%
			\noexpand\@tfor\noexpand\Hy@tempa:=%
			\expandafter\strip@prefix\meaning\Hy@tempa
		}%
		\@tempb\do{%
			\if\Hy@tempa\else
				\if\Hy@tempa\else
					\xdef\newname{\newname\Hy@tempa}%
				\fi
			\fi
		}%
	\endgroup
}%
\newenvironment{enumBib}{%
	\BibTitle
	\advance\@enumdepth \@ne
	\edef\@enumctr{enum\romannumeral\the\@enumdepth}\list
	{\csname biblabel\@enumctr\endcsname}{\usecounter
	{\@enumctr}\def\makelabel##1{\hss\llap{\upshape##1}}}
}{%
	\endlist
}

\def\Chapters#1{\ChapterList#1,LastChapter,}%
\def\LastChapter{LastChapter}%
\def\ChapterList#1,{\def\temp{#1}%
	\ifx\temp\LastChapter
	\else
		\@ifundefined{#1}{%
		}{%
			\def\Semafor{on}
		}
		\expandafter\ChapterList
	\fi
}%
\newcommand{\BiblioItem}[3]
{
	\def\Semafor{off}
	\Chapters{#1}
	\ifx\Semafor\ValueOn
		\ifx\IndexState\ValueOff
			\begin{enumBib}
			\def\IndexState{on}
		\fi
		\item \label{bibitem: #2}#3%
	\fi
}
\newcommand{\OpenBiblio}
{
	\def\IndexState{off}
}
\newcommand{\CloseBiblio}
{
	\ifx\IndexState\ValueOn
		\end{enumBib}
		\def\IndexState{off}
	\fi
}

\def\StartCite{[}%
\def\citeBib#1{\showCiteBib#1,endCite,}%
\def\endCite{endCite}%
\def\showCiteBib#1,{\def\temp{#1}%
\ifx\temp\endCite
]%
\def\StartCite{[}%
\else
	\StartCite\ref{bibitem: #1}%
	\def\StartCite{, }%
\expandafter\showCiteBib%
\fi}%

\makeatother 
\newcommand{\arp}{\ar @{-->}}
\newcommand{\ars}{\ar @{.>}}
\newcommand{\bundle}[4]
{
	\def\tempa{}%
	\def\tempb{#3}%
	\def\tempc{#1}%
	\ifx\tempa\tempb
		\ifx\tempa\tempc
			#2
		\else
			\xymatrix{#2:#1\arp[r]&#4}
		\fi
	\else
		\ifx\tempa\tempc
			#2[#3]
		\else
			\xymatrix{#2[#3]:#1\arp[r]&#4}
		\fi
	\fi
}
\newcommand{\AddIndex}[2]%
{%
	{\bf #1}%
	\label{index: #2}%
}%
\newcommand{\Index}[3]%
{%
	\def\Semafor{off}%
	\Chapters{#1}%
	\ifx\Semafor\ValueOn%
		\def\tempa{}%
		\def\tempb{#3}%
		\ifx\IndexState\ValueOff%
			\begin{theindex}%
			\def\IndexState{on}%
		\fi%
		\ifx\IndexSpace\ValueOn%
			\indexspace%
			\def\IndexSpace{off}%
		\fi%
		\item #2%
		\ifx\tempa\tempb%
		\else%
			\ \pageref{index: #3}%
		\fi%
	\fi%
}%
\newcommand{\SubIndex}[3]
{
	\def\Semafor{off}
	\Chapters{#1}
	\ifx\Semafor\ValueOn
		\subitem #2 \pageref{index: #3}
	\fi
}%

\makeatletter
\newcommand{\Symb}[3]
{
	\def\Semafor{off}
	\Chapters{#1}
	\ifx\Semafor\ValueOn
		\ifx\IndexState\ValueOff
			\begin{theindex}
			\def\IndexState{on}
		\fi
		\ifx\IndexSpace\ValueOn
			\indexspace
			\def\IndexSpace{off}
		\fi
		\item $\@nameuse{ViewSymbol#3}$\ \ #2
		\@nameuse{RefSymbol#3}%
	\fi
}
\def\CiteBibNote{\footnote[0]}
\makeatother

\newcommand{\SetIndexSpace}%
{%
	\def\IndexSpace{on}%
}%
\def\ValueOff{off}
\def\ValueOn{on}

\newcommand{\OpenIndex}
{
	\def\IndexState{off}
}
\newcommand{\CloseIndex}
{
	\ifx\IndexState\ValueOn
		\end{theindex}
		\def\IndexState{off}
	\fi
}
\newcommand\epigraph[2]
{
	\begin{flushright}
	\begin{minipage}{200pt}
		\Memos{#1}\CiteBibNote{#2}
	\end{minipage}
	\end{flushright}
}
\def\Memos#1{\MemoList#1//LastMemo//}%
\def\LastMemo{LastMemo}%
\def\MemoList#1//{\def\temp{#1}%
	\ifx\temp\LastMemo
	\else%
		\par\setlength{\parindent}{12pt}\textcolor{blue}{#1}%
		\expandafter\MemoList%
	\fi%
}%

\def\texPrefaceAlgebra{}
\ifx\PrintBook\Defined
\else
\xRefDef{0701.238}{math/pdf/0701/0701238}
\begin{document}
\title{Fibered $\mathcal{F}$\Hyph Algebra}

\pdfbookmark[1]{Fibered F-Algebra}{TitleEnglish}
\begin{abstract}
The concept of $\mathcal{F}$\Hyph algebra and its representation can be extended to
an	arbitrary bundle. We define operations of fibered $\mathcal{F}$\Hyph algebra
in fiber. The paper presents the representation theory of
of fibered $\mathcal{F}$\Hyph algebra as well as a comparison of representation of
$\mathcal{F}$\Hyph algebra and of representation of fibered $\mathcal{F}$\Hyph algebra.
\end{abstract}

\maketitle
\fi

\def\texPrefaceAlgebra{}
\def\texPreface_Algebra{}
\ifx\texPrefaceRelation\Defined
\ifx\PrintBook\Defined
				\chapter{Preface}
\fi
\epigraph{24 - 25 August 1883
//
Pushkin once said in the circle of his friends: "Imagine what my Tat'iana has done
- she's got married. I should never have expected that of her." I could say just
the same about Anna Karenina. My characters sometimes do things that I would not wish.
In general they do what is ordinarily done in actual life and not what I want.
//
- Reminiscences of V. I. Alekseev. The date is that given by G. A. Rusanov, who
records the same saying. According to Alekseev it refers to Anna's suicide.
Tat'iana is the heroine of Pushkin's "novel in verse," Eugenii Onegin.}
{\citeBib{Tolstoi about Anna Karenina}, p 51}
\fi

\ifx\texPrefaceAlgebra\Defined
Theory of representation of algebra has long and extensive history.
During XX century representation theory became an integral part of different applications.
Transition from algebra to algebra bundle opens new opportunities.
I have ventured to write this
\ifx\PrintBook\Defined
book
\else
paper
\fi
where I want to discover new properties of algebra bundle.
\fi

\ifx\texPrefaceRelation\Defined
This story began many years ago.
When I was 15 years old I got some money for minor expenses.
I started buying math and physics books. When
I entered university I collected library of favorite books.
Among these books there was book "Universal Algebra" by Cohn
(\citeBib{Cohn: Universal Algebra})
that miraculously remained unsold. I was not thinking too hard that
it might be too early for me to read this. Reading this book
I had fallen in love with algebra.
However I dedicated my life and research to geometry 
on the edge of geometry and physics.

After more then 30 years I suddenly returned to this book.
I became curious about
definition of universal algebra in fiber of bundle.
When I started to write paper \citeBib{0702.561} I could swear that
something similar I had read when I was young. I found the book
that, as I supposed, was origin. However this book was dedicated only
to vector bundles. When I started to write this paper I finally
realized that I did not read something similar. I could not leave such facts unnoticed.

It turns out that I thought about this when I was young. However why did I wait so long time?
Why I initiated this research now? I will never get an answer to the
first question. However the answer to the second question is very simple.
For a long time I studied reference frame in general relativity.
My interest was not only classical case, but possible deviations of geometry as well.
Statements obtained in \citeBib{0702.561} show  that
in the field of my research there is a lot of unanswered questions.

I supposed to dedicate this paper to
fibered equivalence relation, because I have certain interest to it
in the future.
However, the need for clear statements involved new definitions.
Then events became unpredictable.
I expected the paper to be extremely concise and
written for a month. However this paper severely takes all my time,
changes the title and direction of the research.\footnote{Thrill
of hunt is one of the strongest passions of mankind.
I catch myself that I keep solving problems which are new for me.}
The paper dedicated to fibered equivalence relation turns into paper
dedicated to fibered binary relations.

Fibered relation is one of the most complicated subjects in the theory of fibered
algebra. Since an operation is a map,
we extend unambiguously the definition of the operation to bundle
and its sections and demand that operation is smooth.
Relation is the subset of Cartesian product.
Assuming that we defined relation only in the fiber, we are losing relationship between fibers.

I decided to repeat the procedure to determine a relation in
universal algebra.
When I started to study fibered correspondence, I realized, that I need
to change definitions in \citeBib{0702.561}.
The definition
\xRef{0702.561}{definition: reduced Cartesian product of bundle}
generates too narrow framework
to define a fibered correspondence.
May be sometimes it is enough to define correspondence only in
fiber\footnote{In particular, we define a fibered relation
introducing relation in a fiber}, however we losing
fibered morphisms.
The analysis of this situation exposes the myth
with which I have comfortably lived for all those years.
More exactly, all this time I have been trying to work out what
the base of fibered map is like. Is it an injection or an arbitrary map?
No one definition
gives clear answer on this question.
For the simplicity of perception I supposed that
base of fibered map is
injection. Actually, since we do not determine type of map of a base,
this map may be arbitrary. This leads to more wide
definition (definition
\xRef{0702.561}{definition: Cartesian product of bundle})
of Cartesian product of bundles.
On the other hand, the definition
\xRef{0702.561}{definition: Cartesian product of bundle}
presents problems to determine fibered algebra. This brings to necessity
to use two definitions of Cartesian product of bundle.
Similar considerations bring to two definitions of fibered correspondence.

В конце концов я вернулся к отношению эквивалентности. Но мой труд
не пропал даром. Мой взгляд на проблему изменился.

From definition \ref{definition: Fibered Correspondence}
it follows that notion of continuity is important in
definition of fibered correspondence. Вслед за этим определением
я разъясняю, что означает непрерывность соответствия.
\fi


\ifx\texPrefaceAlgebra\Defined
Since a cut of the bundle
may be not defined on the whole bundle, all statements assume
a specific domain.
Such statements are modelled on statement
about existance of trivial tangent bundle on manifold.

However there exists other group of statements restricting the
domain of fibered
$\mathcal{F}$\Hyph algebra. I give appropriate examples in the text.
The explanation follows.

We suppose that transormation of fiber caused by parallel transfer
is one-to-one map.
Continuously moving along base,
we continuously move from one fiber to other.
Assumption that map between fibers is homeomorphism
waranties continuous deformation of fiber.

If we choose the point in fiber, then the trajectory
of its movement when its projection moves along
base is parallel to base.
In differential geometry such lines are called horisontal.
We also will use this definition.

Since we assume that there is structure of
$\mathcal{F}$\Hyph algebra on fiber, then we expect that
corresponding map is isomorphism
of $\mathcal{F}$\Hyph algebra.
Continuity allows to save considered structures when we use
parallel transfer, it allows to make smooth transfer from fiber to fiber.

This picture works well in the small. When we consider
finite intervals on base,
continuity becomes responsible for
impossibility to extend the structure of $\mathcal{F}$\Hyph algebra
as far as we please.
For instance, there may apear points where homeomorphism is
broken. It happens when horizontal lines have intersection or
topological propertyes of fiber change. Corresponding fiber is
called degenerate, and its projection is called point of degeneracy.

It is not easy to say how many points of degeneracy there are.
It is clear by intuition that this set is small in comparison with the base.
However, this set may prove to be essential
for the study of geometry of bundles or physical processes
associated with this bundle.

The problem to extend the solution of differential equation
is one of such events in the theory of differential equations.
At the same time, there exist two types of solution of
differential equation.
Regular solution belongs to family of functions
dependent on arbitrary constants.
Singular solution is envelope of family of regular solutions.

The problem of describing fibers of bundle, regardless whether
they are degenerate or not, has an interesting solution. Any path on the base of bundle
is map of interval $I=[0,1]$ into base of bundle.
Let us assume that fibers of bundle are not homeomorphic, but homotopic.

The holonomy group of bundle also has constrains
for structure of fibered $\mathcal F$\Hyph algebra
$\mathcal A=\bundle{}pA{}$. It is natural to assume that using
parallel transfer we have homomorphism of
$\mathcal F$\Hyph algebra from one fiber into another.
Therefore, we assume that transformation caused by
parallel transfer along loop
is homomorfism of $\mathcal F$\Hyph algebra. Thus,
everything is fine
when the holonomy group of bundle $\mathcal A$ is
subgroup of group of homomorphisms of $\mathcal F$\Hyph algebra
$A$. In this case
fibered $\mathcal F$\Hyph algebra $\mathcal A$ is called
holonomic.
Otherwise
fibered $\mathcal F$\Hyph algebra $\mathcal A$ is called
anholonomic.

From theory of vector bundles we know that there
exist fibered $\mathcal F$\Hyph algebra which is
not holonomic. At the same time, the theory of vector bundles
has answer how we can work with
anholonomic fibered $\mathcal F$\Hyph algebra.

We apply this remark also to the theory of representations of
fibered $\mathcal F$\Hyph algebra.

A new design is illustrated through the use of corresponding diagrams.
\fi

Where it is possible
I use the same notation for operations and relations
as we use them in the set theory. It does not bring to
ambiguity because we use different notation for set and
bundle. I use the same letter in different alphabets to denote
bundle and fiber.

We assume that projection of bundle, section and fibered map are smooth maps.

\def\texCartesian{}
\ifx\PrintBook\Defined
				\chapter{Cartesian Product of Bundles}
\fi

			\section{Bundle}

Let $M$ be a manifold and
\symb{\bundle{\mathcal{E}}pEM}0{bundle}
	\begin{equation}
	\label{eq: bundle, definition}
\ShowSymbol{bundle}
	\end{equation}
be a bundle over $M$ with fiber $E$.\footnote{Since I have
deal with different bundles I follow next agreement.
I use the same letter in different alphabets for notation of bundle and fiber.}
The symbol $\bundle{}pE{}$ means that
$E$ is a typical fiber of the bundle.
Set $\mathcal{E}$ is domain of
map $\bundle{}pE{}$.
Set $M$ is range of
map $\bundle{}pE{}$.
We identify the smooth map $\bundle{}pE{}$
and the bundle \eqref{eq: bundle, definition}.
Denote by
\symb{\Gamma(\bundle{}pE{})}1{set of sections of bundle}
the set of sections of bundle $\bundle{}pE{}$.

\AddIndex{Cartesian power $A$ of set $B$}{Cartesian power of set}
is the set \symb{B^A}1{Cartesian power of set}
of mappings $f:A\rightarrow B$
(\citeBib{Cohn: Universal Algebra}, page 5).
Let us consider subsets of $B^A$ of the form
\[
W_{K,U}=\{f:A\rightarrow B|f(K)\subset U\}
\]
where $K$ is compact subset of space $A$,
$U$ is open subset of space $B$.
Sets $W_{K,U}$ form base of topology
on space $B^A$. This topology is called 
\AddIndex{compact\hyph open topology}{compact open topology}.
Cartesian power $A$ of set $B$
equipped by compact\hyph open topology is called
\AddIndex{mapping space}{mapping space}
(\citeBib{Maunder: Algebraic Topology}, page 213).

According to \citeBib{Maunder: Algebraic Topology}, стр. 214,
given spaces $A$, $B$, $C$, $D$ and
mappings $f:A\rightarrow C$, $g:D\rightarrow B$
we define morphism of mapping spaces
\[
g^f:D^C\rightarrow B^A
\]
by law
\begin{align*}
g^f(h)&=fhg&h:&C\rightarrow D&g^f(h):&A\rightarrow B
\end{align*}
Thus, we can represent the morphism of mapping spaces
using diagram
	\[
\xymatrix{
A\ar[rr]^f\ar[d]_{g^f(h)}&&C\ar[d]^h\\
B&&D\ar[ll]^g
}
	\]

Set $\Gamma(\mathcal{E})$ is subset
of set $\mathcal{E}^M$.
This is why for set of sections we can use definitions established
for mapping set.
This is why for set of sections we can use methods
defined for  mapping space.
We define sets $W_{K,U}$ by law
\[
W_{K,U}=\{f\in\Gamma(\mathcal{E})|f(K)\subset U\}
\]
where $K$ is compact subset of space $M$,
$U$ is open subset of space $\mathcal{E}$.

	\begin{remark}
I use arrow $\xymatrix{\arp[r]&}$ to represent
projection of bundle on diagram.
	\qed
	\end{remark}

	\begin{remark}
I use arrow $\xymatrix{\ars[r]&}$ to represent
section of bundle on diagram.
	\qed
	\end{remark}

Let $f$ be a fibered map
from $\mathcal{E}$ to $\mathcal{E}'$
	\[
\xymatrix{
\mathcal{E}\arp[d]^{\bundle{}pE{}}\ar[rr]^f
&&\mathcal{E}'\arp[d]^{\bundle{}{p'}{E'}{}}\\
M\ar[rr]^F&&M'
}
	\]
The map $F$ is the \AddIndex{base of map}{base of map}
$f$.
The map $f$ is the \AddIndex{lift of map}{lift of map} $F$.

Suppose map $F$ is bijection.
Then the map $f$ defines morphism $f^{F^{-1}}$ of spaces of sections
$\Gamma(\bundle{}pE{}))$ to $\Gamma(\bundle{}{p'}{E'}{}))$
	\[
\xymatrix{
\mathcal{E}\ar[rr]^f&&\mathcal{E}'
&u'=f^{F^{-1}}(u)=fuF^{-1}\\
M\ars[u]^u="U"\ar[rr]^F&&M'\ars[u]_{u'}="UP"
\ar @{=>}^{f^{F^{-1}}} "U";"UP"
}
	\]%
It is enough to prove continuity of $F^{-1}$
to prove continuity of $u'$. However this is evident, because $F$
is continuous bijection.

Since $F=\mathrm{id}$, then $\mathrm{id}^{-1}=\mathrm{id}$.
In this case
we use notation $f^{\mathrm{id}}$ for morphism of spaces of sections.
It is evident, that
\[
f^{\mathrm{id}}(u)=fu
\]

\ifx\texFuture\Defined
Прежде всего возникает вопрос. Что мы называем расслоением? Когда
мы говорим о расслоении, мы подразумеваем непрерывность базы и
непрерывность слоя. Если слой имеет дискретную топологию, то мы
называем соответствующее множество не расслоением, а слоением.
Расслоения и слоения очень сильно отличаются, и потому мы
используем различные методы для их изучения.
 
Тем не менее, если отвлечься от топологии слоя, эти объекты имееют
нечто общее. Кроме того, нужно иметь в виду следующее.
\begin{itemize}
\item Алгебра может быть определена не только на непрерывном множестве.
Конечные и дискретные группы, например, важны при изучении
многих геометрических объектов и физических явлений.
\item Фактор множество топологического пространства может оказаться
конечным множеством, тем самым разрушая топологию.
В то же время совсем не очевидно, почему мы не должны
рассматривать подобного рода отношение эквивалентности.
\end{itemize}

На данном этапе для меня топология слоя несущественна, но важна
непрерывность базы. Принимая во внимание высказанные выше соображения,
я не буду делать различия между расслоением и накрытием
и буду называть и то, и другое расслоением.
\fi

			\section{Cartesian Product of Bundles}
			\label{section: Cartesian Product of Bundдles}
		\begin{remark}
		\label{remark: map to Cartesian product}
Let $A_1\times...\times A_n$ be Cartesian product
of sets $A_1$, ..., $A_n$
	\[
\xymatrix{
&A_1\times...\times A_n\ar[rd]^{p_n}\ar[ld]_{p_1}\ar[d]&\\
A_1&...&A_n
}
	\]
According to \citeBib{Serge Lang} we can represent the map
\[
f:A\rightarrow A_1\times...\times A_n
\]
as tuple \symb{f=(f_1,...,f_n)}1{map to Cartesian product}
where $f_i=fp_i$
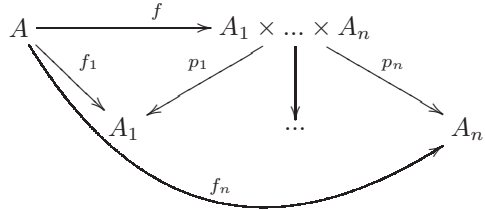
	\[
\xymatrix{
A\ar[rr]^f\ar[rd]^{f_1}\ar@/_4pc/[rrrd]^{f_n}
&&A_1\times...\times A_n\ar[rd]^{p_n}\ar[ld]_{p_1}\ar[d]&\\
&A_1&...&A_n
}
	\]
	 \qed
		\end{remark}

Let $\bundle{\mathcal{E}_i}{p_i}{E_i}{M_i}$,
$i=1$, ..., $n$ be the set of bundles.
For any $i$ let $\{U_{i\alpha_i}\}$ be a cover of
$M_i$ such that for any $U_{i\alpha_i}$ there exists
a local chart $\varphi_{i\alpha_i}$
of the bundle $\bundle{}{p_i}{E_i}{}$
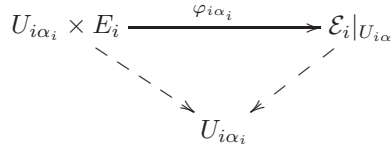
	\[
\xymatrix{
U_{i\alpha_i}\times E_i\ar[rr]^{\varphi_{i\alpha_i}}
\arp[rd]&&\mathcal{E}_i|_{U_{i\alpha}}\arp[ld]\\
&U_{i\alpha_i}
}
	\]
	 \begin{definition}
	 \label{definition: Cartesian product of bundle}
For any point $x_i\in M_i$ there exist
open set $U_{i\alpha_i}$, $i=1$, ..., $n$ 
such, that $x_i\in U_{i\alpha_i}$.
By $(x_i,a_i)$, $x_i\in U_{i\alpha_i}$, $a_i\in E_i$
denote a point of set
$\mathcal{E}_i|_{U_{i\alpha_i}}$.
For tuple $\alpha=(\alpha_1,...,\alpha_n)$ we introduce
trivial bundle
\[\bundle{\mathcal{E}_\alpha}{e_\alpha}{}{U_\alpha=U_{1\alpha_1}\times...\times U_{n\alpha_n}}\]
with $\mathcal{E}_\alpha$ consisting of tuples $(x_1,...,x_n,a_1,...,a_n)$.
Cartesian product $\prod_{i=1}^n E_i$
is a fiber of this bundle.
Cartesian product $\prod_{i=1}^n U_{i\alpha_i}$
is a base of this bundle.
	 \qed
	 \end{definition}

Continuity of projection of this bundle follows from
Corollary 1 of Proposition 1
(\citeBib{Bourbaki: General Topology 1}, page 44).

To define a bundle over manifold $\prod_{i=1}^nM_i$,
we need to define gluing functions. Let $\psi_{i\alpha_i\beta_i}$
be gluing functions of bundle $\bundle{}{p_i}{E_i}{}$ 
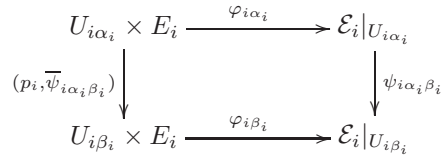
	\[
\xymatrix{
U_{i\alpha_i}\times E_i\ar[rr]^{\varphi_{i\alpha_i}}
\ar[d]_{(p_i,\overline{\psi}_{i\alpha_i\beta_i})}
&&\mathcal{E}_i|_{U_{i\alpha_i}}
\ar[d]^{\psi_{i\alpha_i\beta_i}}\\
U_{i\beta_i}\times E_i\ar[rr]^{\varphi_{i\beta_i}}
&&\mathcal{E}_i|_{U_{i\beta_i}}
}
	\]
Suppose $U_\alpha\cap U_\beta\ne\emptyset$. Then we
define gluing function using diagram
	
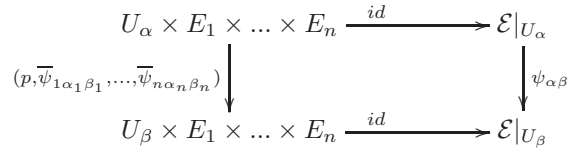
\[
\xymatrix{
U_\alpha\times E_1\times...\times E_n
\ar[rr]^{id}
\ar[d]_{(p,\overline{\psi}_{1\alpha_1\beta_1},...,
\overline{\psi}_{n\alpha_n\beta_n})}
&&\mathcal{E}|_{U_\alpha}\ar[d]^{\psi_{\alpha\beta}}\\
U_\beta\times E_1\times...\times E_n
\ar[rr]^{id}
&&\mathcal{E}|_{U_\beta}
}
	\]

Bundle
\symb{\bundle{}{p_1}{E_1}{}\times...\times\bundle{}{p_n}{E_n}{}}
0{Cartesian product of bundles, definition 1}
\symb{\bundle{}{\prod_{i=1}^n{p_i}}{E_i}{}}
0{Cartesian product of bundles, definition 2}
\[
\ShowSymbol{Cartesian product of bundles, definition 1}=
\ShowSymbol{Cartesian product of bundles, definition 2}=
\bundle{\mathcal{E}}p{\prod_{i=1}^nE_i}M
\]
is called \AddIndex{Cartesian product of bundles}{Cartesian product of bundles}
$\bundle{}{p_i}{E_i}{}$.
We also speak that the total space $\mathcal{E}$
is \AddIndex{Cartesian product of total spaces}{Cartesian product of total spaces}
$\mathcal{E}_i$ and use notation
\symb{\prod_{i=1}^n{\mathcal{E}_i}}
0{Cartesian product of total spaces, definition 2}
\symb{\mathcal{E}_1\times...\times\mathcal{E}_n}
0{Cartesian product of total spaces, definition 1}
\[
\mathcal{E}=
\ShowSymbol{Cartesian product of total spaces, definition 1}=
\ShowSymbol{Cartesian product of total spaces, definition 2}
\]

		\begin{remark}
		\label{remark: Cartesian product, section}
According to remark \ref{remark: map to Cartesian product}
we can represent a section of Cartesian product of bundles
$\bundle{}{p_1}{E_1}{}\times...\times\bundle{}{p_n}{E_n}{}$
as tuple of sections $a=(a_1,...,a_n)$.
		\qed
		\end{remark}

We will use following diagrams to represent Cartesian product of bundles
	
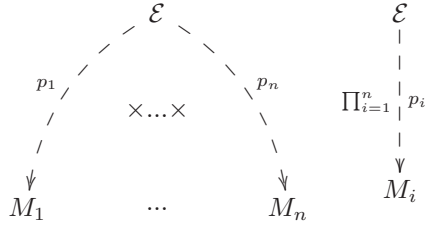
\[
\xymatrix{
&\mathcal{E}\arp@/^1pc/[rdd]^{p_n}\arp@/_1pc/[ldd]_{p_1}&\\
&\times...\times&\\
M_1&...&M_n
}
\ \ \xymatrix{
\mathcal{E}\arp[dd]^{p_i}_{\prod_{i=1}^n}\\
\\
M_i
}
	\]%
On the diagram, the arrows connected by either symbol $\times$ or $\prod$
denote the arrow corresponding to projection of bundle
$\mathcal{E}$. The notation is intended to show the structure of map.

Product of any two fibers is defined in
the Cartesian product of bundles. As we see below, such bundle reproduces the structure of the base.
This substantially
restricts the application of Cartesian product of bundles.

			\section{Reduced Cartesian Product of Bundles}
			\label{section: reduced Cartesian Product of Bundles}

Since bundles are defined over common base we can change definition of Cartesian product of bundles.

Let $\bundle{\mathcal{E}_i}{p_i}{E_i}M$,
$i=1$, ..., $n$ be the set of bundles over manifold $M$.
For any $i$ let $\{U_{i\alpha_i}\}$ be a cover of
$M$ such that for any $U_{i\alpha_i}$ there exists
a local chart $\varphi_{i\alpha_i}$
of the bundle $\bundle{}{p_i}{E_i}{}$
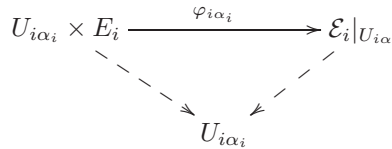
	\[
\xymatrix{
U_{i\alpha_i}\times E_i\ar[rr]^{\varphi_{i\alpha_i}}
\arp[rd]&&\mathcal{E}_i|_{U_{i\alpha}}\arp[ld]\\
&U_{i\alpha_i}
}
	\]
	 \begin{definition}
	 \label{definition: reduced Cartesian product of bundle}
For any point $x\in M$ there exist
open sets $U_{i\alpha_i}$, $i=1$, ..., $n$ 
such, that $x\in U_{i\alpha_i}$.
By $(x,a_i)$, $x\in U_{i\alpha_i}$, $a_i\in E_i$
denote a point of set
$\mathcal{E}_i|_{U_{i\alpha_i}}$.
For tuple $\alpha=(\alpha_1,...,\alpha_n)$ we introduce
trivial bundle
$\bundle{\mathcal{E}_\alpha}{p_\alpha}{}{U_\alpha=\bigcap_{i=1}^nU_{i\alpha_i}}$
with $\mathcal{E}_\alpha$ consisting of tuples $(x,a_1,...,a_n)$.
Cartesian product $\prod_{i=1}^n E_i$
is a fiber of this bundle.
The set $U_\alpha$ is a base of this bundle.
	 \qed
	 \end{definition}

According to definition, we can represent
the bundle over set $U_\alpha$ as
\[
U_\alpha\times\prod_{i=1}^n E_i
\]
According to \citeBib{Bourbaki: General Topology 1}, page 44,
given $U$ belongs to the base of topology of space $U_\alpha$,
then 
$U\times\prod_{i=1}^n E_i$
belongs to the base of topology of space $\mathcal{E}_\alpha$.
We can represent set $U$ as $U=\cup_{i\in I}U_i$,
where $U_i$ belong to the base of topology of space $M$,
if $U$ is open set of space $M$.
Accordingly, the set
\[
p_\alpha^{-1}(U)=U\times\prod_{i=1}^n E_i
\]
can be represented as
\[
p_\alpha^{-1}(U)=\bigcup_{i\in I}(U_i\times\prod_{i=1}^n E_i)
\]
and it is an open set.
Therefore, projection $p_\alpha$ is continuous mapping.

To define a bundle over manifold $M$,
we need to define gluing functions. Let $\psi_{i\alpha_i\beta_i}$
be gluing functions of bundle $\bundle{}{p_i}{E_i}{}$ 
	\[
\xymatrix{
U_{i\alpha_i}\times E_i\ar[rr]^{\varphi_{i\alpha_i}}
\ar[d]_{(p_i,\overline{\psi}_{i\alpha_i\beta_i})}
&&\mathcal{E}_i|_{U_{i\alpha_i}}
\ar[d]^{\psi_{i\alpha_i\beta_i}}\\
U_{i\beta_i}\times E_i\ar[rr]^{\varphi_{i\beta_i}}
&&\mathcal{E}_i|_{U_{i\beta_i}}
}
	\]
Suppose $U_\alpha\cap U_\beta\ne\emptyset$. Then we
define gluing function using diagram
	\[
\xymatrix{
U_\alpha\times E_1\times...\times E_n
\ar[rr]^{id}
\ar[d]_{(p,\overline{\psi}_{1\alpha_1\beta_1},...,
\overline{\psi}_{n\alpha_n\beta_n})}
&&\mathcal{E}|_{U_\alpha}\ar[d]^{\psi_{\alpha\beta}}\\
U_\beta\times E_1\times...\times E_n
\ar[rr]^{id}
&&\mathcal{E}|_{U_\beta}
}
	\]

Bundle
\symb{\bundle{}{p_1}{E_1}{}\odot...\odot\bundle{}{p_n}{E_n}{}}
0{reduced Cartesian product of bundles, definition 1}
\symb{\bundle{}{\bigodot_{i=1}^n{p_i}}{E_i}{}}
0{reduced Cartesian product of bundles, definition 2}
\[
\ShowSymbol{reduced Cartesian product of bundles, definition 1}=
\ShowSymbol{reduced Cartesian product of bundles, definition 2}=
\bundle{\mathcal{E}}p{\prod_{i=1}^nE_i}M
\]
is called \AddIndex{reduced Cartesian product of bundles}{reduced Cartesian product of bundles}
$\bundle{}{p_i}{E_i}{}$.
We also speak that the total space $\mathcal{E}$
is \AddIndex{reduced Cartesian product of total spaces}{reduced Cartesian product of total spaces}
$\mathcal{E}_i$ and use notation
\symb{\bigodot_{i=1}^n{\mathcal{E}_i}}
0{reduced Cartesian product of total spaces, definition 2}
\symb{\mathcal{E}_1\odot...\odot\mathcal{E}_n}
0{reduced Cartesian product of total spaces, definition 1}
\[
\mathcal{E}=
\ShowSymbol{reduced Cartesian product of total spaces, definition 1}=
\ShowSymbol{reduced Cartesian product of total spaces, definition 2}
\]

		\begin{remark}
		\label{remark: reduced Cartesian product, section}
According to remark \ref{remark: map to Cartesian product}
we can represent a section of reduced Cartesian product of bundles
$\bundle{}{p_1}{E_1}{}\odot...\odot\bundle{}{p_n}{E_n}{}$
as tuple of sections $a=(a_1,...,a_n)$.
		\qed
		\end{remark}

We will use following diagrams to represent reduced Cartesian product of bundles
	
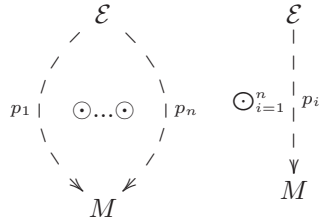
\[
\xymatrix{
\mathcal{E}\arp@/^2pc/[dd]^{p_n}\arp@/_2pc/[dd]_{p_1}\\
\odot...\odot\\
M
}
\ \ \xymatrix{
\mathcal{E}\arp[dd]^{p_i}_{\bigodot_{i=1}^n}\\
\\
M
}
	\]%
On the diagram, the arrows connected by symbol $\odot$
denote the arrow corresponding to projection of bundle
$\mathcal{E}$. The notation is intended to show the structure of map.

In reduced Cartesian product we define product of fibers over selected point.
This makes the structure of product more rich.

	 \begin{definition}
	 \label{definition: Cartesian power of bundle}
For $n\ge 0$ we define
\AddIndex{Cartesian power of bundle}{Cartesian power of bundle}\footnote{Since I use
definition of Cartesian power of bundle only in frame of reduced Cartesian product, I do noot use
respective adjective for power. I use this remark for all following
definitions related to Cartesian power of bundle.}
\symb{\bundle{}{p}{E}{}^n}0{Cartesian power of bundle}
\symb{\mathcal{E}^n}0{Cartesian power of total spaces}
\[
\left\{
\begin{array}{lr}
\bundle{}{p}{E}{}^0
=\bundle{M}{id}{}M
& n=0\\
\ShowSymbol{Cartesian power of bundle}
=\bundle{\ShowSymbol{Cartesian power of total spaces}}
{\bigodot_{i=1}^n\bundle{}{p}{E}{}}{}M
&n>0
\end{array}
\right.
\]
	 \qed
	 \end{definition}

\def\texFiberedAlgebra{}
\ifx\PrintBook\Defined
				\chapter{Fibered \texorpdfstring{$\mathcal{F}$}{F}\Hyph Algebra}
\fi

			\section{Fibered \texorpdfstring{$\mathcal{F}$}{F}\Hyph Algebra}

	 \begin{definition}
An n-ary
\AddIndex{operation on bundle}{operation on bundle}
$\bundle{}pE{}$
is a fibered map
\[
f:\mathcal{E}^n\rightarrow \mathcal{E}
\]
$n$ is \AddIndex{arity of operation}{arity of operation}.
$0$-arity operation is a section of $\mathcal{E}$.
	 \qed
	 \end{definition}

We can represent the operation using the diagram
	
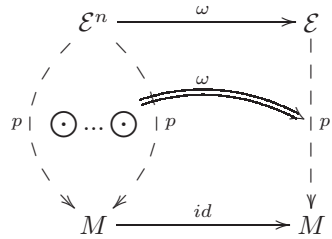
\[
\xymatrix{
\mathcal{E}^n\ar[rr]^\omega\arp@/^2pc/[dd]^p\arp@/_2pc/[dd]_p&&\mathcal{E}\arp[dd]^p\\
\bigodot...\bigodot\ar@{=>}@/^1pc/[rr]^\omega &&\\
M\ar[rr]^{id}&&M
}
	\]

		\begin{theorem}
		\label{theorem: continuity of operation}
Let $U$ be an open set of base $M$.
Suppose there exist trivialization of bundle $\bundle{}pE{}$ over $U$.
Let $x\in M$.
Let $\omega$ be $n$\Hyph ary operation on bundle $\bundle{}pE{}$ and
\[
\omega(p_1,...,p_n)=p
\]
in the fiber $E_x$.
Then there exist open sets $V\subseteq U$,
$W\subseteq E$, $W_1\subseteq E_1$, ..., $W_n\subseteq E_n$
such, that $x\in V$, $p\in W$, $p_1\in W_1$, ..., $p_n\in W_n$,
and for any $x'\in V$, $p'\in W\cap\omega V$ there exist
$p'_1\in W_1$, ..., $p'_n\in W_n$ such, that
\[
\omega(p'_1,...,p'_n)=p'
\]
in the fiber $E_{x'}$.
		\end{theorem}
		\begin{proof}
According to \citeBib{Bourbaki: General Topology 1}, page 44,
since $V$ belongs to the base
of topology of space $U$ and $W$ belongs to the base
of topology of space $E$, then
set $V\times W$ 
belongs the base of topology of space $\mathcal{E}$.
Similarly, since 
$V$ belongs to the base
of topology of space $U$ and
$W_1$, ..., $W_n$ belong to the base of topology of space $E$,
set $V\times W_1\times ...\times W_n$
belongs the base of topology of space $\mathcal{E}^n$.

Since mapping $\omega$ is continuous, then for an open set
$V\times W$ there exists an open set $S\subseteq \mathcal{E}^n$ such,
that $\omega S\subseteq V\times W$.
Suppose $x'\in V$.
Let
$(x',p')\in\omega S$ be an arbitrary point.
Then there exist such $p'_1\in E_{x'}$, ..., $p'_n\in E_{x'}$,
that
\[
\omega(p'_1,...,p'_n)=p'
\]
in fiber $E_{x'}$.
According to this there exist sets $R$, $R'$
from base of topology of space $U$,
and sets $T_1$, ..., $T_n$,
$T'_1$, ..., $T'_n$ from base of topology of space $E$,
such that $x\in R$, $x'\in R'$,
$p_1\in T_1$, $p'_1\in T'_1$, ..., $p_n\in T_n$, $p'_n\in T'_n$,
$R\times T_1\times ...\times T_n\subseteq S$,
$R'\times T'_1\times ...\times T'_n\subseteq S$.
We proved the theorem since
$W_1=T_1\cup T'_1$, ..., $W_n=T_n\cup T'_n$
are open sets.
%
		\end{proof}

Theorem \ref{theorem: continuity of operation} tells
about continuity of operation $\omega$, however this
theorem tells nothing regarding sets
$W_1$, ..., $W_n$. In particular, it is possible that these sets
are not connected.

We suppose $W=\{p\}$, $W_1=\{p_1\}$, ..., $W_n=\{p_n\}$,
if topology on fiber $A$ is discrete.
This leads one to assume that  in the neighborhood $V$
the operation does not depend on a fiber.
We call the operation $\omega$ locally constant.
However, it is possible that a condition of constancy
is broken on bundle in general.
Thus  the covering space $R\rightarrow S^1$ of the circle $S^1$
defined by $p(t)=(\sin t,\cos t)$ for any $t\in R$ is
bundle over circle with fiber of group of integers.

Let us consider the alternative point of view on the
continuity of operation $\omega$
to get a better understanding of role of continuity
Let us consider the continuity of operation $\omega$
to better see what does it mean. We need to consider
sections, if we want to show that
infinitesimal change of operand when moving along base
causes infinitesimal change of operation.
This change is legal,
because we defined operation on bundle in fiber.

		\begin{theorem}
		\label{theorem: operation over sections}
An n-ary operation on bundle maps sections into section.
		\end{theorem}
		\begin{proof}
Suppose $f_1$, ..., $f_n$ are sections and we define map
	\begin{equation}
f=\omega^{\mathrm{id}}(f_1,...,f_n):M\rightarrow\mathcal{E}
	\label{eq: operation over sections}
	\end{equation}
as
	\begin{equation}
f(x)=\omega(f_1(x),...,f_n(x))
	\label{eq: operation over sections, x}
	\end{equation}
Let $x\in M$ and $u=f(x)$. Let $U$ be a neighborhood of the point $u$ in the range
of the map $f$.

Since $\omega$ is smooth map,
then according to \citeBib{Bourbaki: General Topology 1}, page 44,
for any $i$,
$1\le i\le n$ the set $U_i$ is defined in the range of section $f_i$
such, that $\prod_{i=1}^nU_i$
is open in the range of section $(f_1,...,f_n)$ of the bundle $\mathcal{E}^n$ and
\[
\omega(\prod_{i=1}^nU_i)\subseteq U
\]

Let $u'\in U$.
Since $f$ is a map, then there exist $x'\in M$ such that $f(x')=u'$.
From equation \eqref{eq: operation over sections, x} it follows
that there exist $u'_i\in U_i$, $p(u'_i)=x'$ such, that
$\omega(u'_1,...,u'_n)=u'$.
Since $f_i$ is a section, then there exist a set $V_i\subseteq M$ such, 
that $f_i(V_i)\subseteq U_i$
and $x\in V_i$, $x'\in V_i$. Therefore, the set
\[
V=\cap_{i=1}^nV_i
\]
is not empty, it is open in $M$ and $x\in V$, $x'\in V$.
Thus the map $f$ is smooth and $f$ is the section.
		\end{proof}

We can represent the operation using the diagram
	
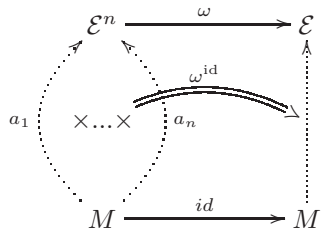
\[
\xymatrix{
\mathcal{E}^n\ar[rr]^\omega&&\mathcal{E}\\
\times...\times\ar@{=>}@/^1pc/[rr]^{\omega^{\mathrm{id}}} &&\\
M\ars@/^2pc/[uu]^{a_1}\ars@/_2pc/[uu]_{a_n}\ar[rr]^{id}&&M\ars[uu]
}
	\]

		\begin{theorem}
		\label{theorem: continuous operation over sections}
$\omega^{\mathrm{id}}$ is continuous on $\Gamma(\mathcal{E})$.
		\end{theorem}
		\begin{proof}
Let us consider a set $W_{K,U}$,
where $K$ is compact set of space $M$,
$U$ is open set of space $\mathcal{E}$.
We can represent set $U$ as $V\times E$,
where $V$ is open set of space $M$, $K\subset V$.
$\omega^{-1}(V\times E)=V\times E^n$ is open set.
Therefore,
\begin{equation}
(\omega^{\mathrm{id}})^{-1}W_{K,V\times E}=W_{K,V\times E^n}
\label{eq: continuous operation over sections}
\end{equation}
From \eqref{eq: continuous operation over sections}
continuity of $\omega^{\mathrm{id}}$ follows.
\end{proof}

	 \begin{definition}
Let $A$ be $\mathcal{F}$\Hyph algebra
(\citeBib{Burris Sankappanavar}).
We can extend $\mathcal{F}$\Hyph algebraic structure from fiber
$A$ to bundle
$\bundle{\mathcal{A}}pAM$.
If operation $\omega$ is defined on $\mathcal{F}$\Hyph algebra $A$
\[
a=\omega(a_1,...,a_n)
\]
then operation $\omega$ is defined on bundle
\[
a(x)=\omega(a_1,...,a_n)(x)=\omega(a_1(x),...,a_n(x))
\]
We say that $\bundle{}pA{}$ is
a \AddIndex{fibered $\mathcal{F}$\Hyph algebra}{fibered F-algebra}.
	 \qed
	 \end{definition}

Depending on the structure we talk for instance about
\AddIndex{fibered group}{fibered group},
\AddIndex{fibered ring}{fibered ring},
or \AddIndex{vector bundle}{vector bundle}.

Main properties of $\mathcal{F}$\Hyph algebra
hold for fibered $\mathcal{F}$\Hyph algebra as well. Proving
appropriate theorems we can refer on this statement.
However in certain cases the proof itself may be of deep interest,
allowing a better view
of the structure of the fibered $\mathcal{F}$\Hyph algebra.
However properties of $\mathcal{F}$\Hyph algebra on the set of sections
are different from properties of $\mathcal{F}$\Hyph algebra in fiber.
For instance, if the product in fiber has inverse element,
it does not mean that the product of sections has inverse element.
Therefore, fibered continuous field generates
ring on the set of sections.
This is the advantage when we consider fibered algebra.
I want also to stress that the operation on bundle is not defined for
elements from different fibers.

Let transition functions $g_{\epsilon\delta}$ determine bundle $\mathcal{B}$
over base $N$.
Let us consider maps $V_\epsilon\in N$ and $V_\delta\in N$,
$V_\epsilon\cap V_\delta\ne\emptyset$.
Point $q\in\mathcal{B}$ has representation $(y,q_\epsilon)$
in map $V_\epsilon$ and representation
$(y,q_\delta)$ in map $V_\delta$.
Therefore,
\[
p_\alpha=f_{\alpha\beta}(p_\beta)
\]
\[
q_\epsilon=g_{\epsilon\delta}(q_\delta)
\]
When we move from map $U_\alpha$ to map $U_\beta$
and from map $V_\epsilon$ to map $V_\delta$,
representation of correspondence changes according to the law
\[
(x,y,p_\alpha,q_\epsilon)=(x,y,f_{\alpha\beta}(p_\beta),g_{\epsilon\delta}(q_\delta))
\]
This is consistent with the transformation
when we move from map $U_\alpha\times V_\epsilon$ to map $U_\beta\times V_\delta$
in the bundle $\mathcal{A}\times\mathcal{B}$.

		\begin{theorem}
		\label{theorem: operation in two neighborhood}
Let transition functions $f_{\alpha\beta}$ determine
fibered $\mathcal{F}$\Hyph algebra $\bundle{\mathcal{A}}pAM$
over base $M$.
Then transition functions $f_{\alpha\beta}$ are
homomorphisms of $\mathcal{F}$\Hyph algebra $A$.
		\end{theorem}
		\begin{proof}
Let $U_\alpha\in M$ and $U_\beta\in M$,
$U_\alpha\cap U_\beta\ne\emptyset$ be
neighborhoods where fibered $\mathcal{F}$\Hyph algebra $\bundle{}pA{}$ is trivial.
Let
	\begin{equation}
a_\beta=f_{\beta\alpha}(a_\alpha)
	\label{eq: operation in two neighborhood}
	\end{equation}
be map from bundle $\bundle{}pA{}|_{U_\alpha}$ into
bundle $\bundle{}pA{}|_{U_\beta}$.
Let $\omega$ be $n$\hyph ary operation and points $e_1$, ..., $e_n$
belong to fiber $A_x$, $x\in U_1\cap U_2$. Suppose
	\begin{equation}
e=\omega(e_1,...,e_n)
	\label{eq: operation in two neighborhood, expression}
	\end{equation}

We represent point $e\in\bundle{}pA{}|_{U_\alpha}$ as
$(x,e_\alpha)$ and point $e_i\bundle{}pA{}|_{U_\alpha}$ as
$(x,e_{i\alpha})$.
We represent point $e\in\bundle{}pA{}|_{U_\beta}$ as
$(x,e_\beta)$ and point $e_i\in\bundle{}pA{}|_{U_\beta}$ as
$(x,e_{i\beta})$. According to \eqref{eq: operation in two neighborhood}
	\begin{equation}
e_\beta=f_{\beta\alpha}(e_\alpha)
	\label{eq: operation in two neighborhood, e}
	\end{equation}
	\begin{equation}
e_{i\beta}=f_{\beta\alpha}(e_{i\alpha})
	\label{eq: operation in two neighborhood, ei}
	\end{equation}
According to \eqref{eq: operation in two neighborhood, expression},
the operation in the bundle $A_x$ over neihgborhood $U_\beta$ is
	\begin{equation}
e_\beta=\omega(e_{1\beta},...,e_{n\beta})
	\label{eq: operation in two neighborhood, fiber}
	\end{equation}
Substituting \eqref{eq: operation in two neighborhood, e},
\eqref{eq: operation in two neighborhood, ei} into
\eqref{eq: operation in two neighborhood, fiber} we get
\[
f_{\beta\alpha}(e_\alpha)=\omega(f_{\beta\alpha}(e_{1\alpha}),...,f_{\beta\alpha}(e_{n\alpha}))
\]
This proves that $f_{\beta\alpha}$ is homomorphism of $\mathcal{F}$\Hyph algebra.
		\end{proof}

	\begin{definition}
	\label{definition: homomorphism of fibered F-algebras} 
Let $\bundle{\mathcal{A}}{p}{A}M$
and $\bundle{\mathcal{A}'}{p'}{A'}{M'}$
be two fibered $\mathcal{F}$\Hyph algebras. Bundle map
\symb{f}0{homomorphism of fibered F-algebras}
\[
\ShowSymbol{homomorphism of fibered F-algebras}:\mathcal{A}\rightarrow\mathcal{A}'
\]
is called
\AddIndex{homomorphism of fibered $\mathcal{F}$\Hyph algebra}{homomorphism of fibered F-algebras}
if respective fiber map
\symb{f_x}0{fiber map}
\[
\ShowSymbol{fiber map}:A_x\rightarrow A_{x'}'
\]
is homomorphism of $\mathcal{F}$\Hyph algebra $A$.
	 \qed
	 \end{definition}

	\begin{definition}
	\label{definition: isomorphism of fibered F-algebras} 
Let $\bundle{\mathcal{A}}{p}{A}M$
and $\bundle{\mathcal{A}'}{p'}{A'}{M'}$
be two fibered $\mathcal{F}$\Hyph algebras. Homomorphism of fibered $\mathcal{F}$\Hyph algebras
$f$
is called
\AddIndex{isomorphism of fibered $\mathcal{F}$\Hyph algebras}{isomorphism of fibered F-algebras}
if respective fiber map
\symb{f_x}0{fiber map}
\[
\ShowSymbol{fiber map}:A_x\rightarrow A_{x'}'
\]
is isomorphism of $\mathcal{F}$\Hyph algebra $A$.
	 \qed
	 \end{definition}

	\begin{definition}
Let $\bundle{\mathcal{A}}pAM$ be
an $\mathcal{F}$\Hyph fibered F-algebra and $A'$ be $\mathcal{F}$\Hyph subalgebra
of the $\mathcal{F}$\Hyph algebra $A$.
An fibered $\mathcal{F}$\Hyph algebra
$\bundle{\mathcal{A}'}p{A'}M$
is a \AddIndex{fibered $\mathcal{F}$\Hyph subalgebra}{fibered F-subalgebra}
of the fibered $\mathcal{F}$\Hyph algebra
$\bundle{}pA{}$ if
homomorphism of fibered $\mathcal{F}$\Hyph algebras
$\mathcal{A}'\rightarrow\mathcal{A}$
is fiber embedding.
	 \qed
	 \end{definition}

The homomorphism of fibered $\mathcal{F}$\Hyph algebra is essential part of this definition.
We can break continuity,
if we just limit ourselves to the fact of the existence
of $\mathcal{F}$\Hyph subalgebra in each fiber.

We defined an operation based reduced Cartesian product of bundles.
Suppose we defined an operation based Cartesian product of bundles.
Then the operation is defined for any elements of the bundle. However,
since $p(a_i)=p(b_i)$, $i=1$, ..., $n$, then $p(\omega(a_1,...,a_n))=p(\omega(b_1,...,b_n))$.
Therefore, the operation is defined between fibers. We can map this operation to base
using projection. This structure is not different from quotient $\mathcal{F}$\Hyph algebra
and does not create new element in bundle theory. The same time 
mapping between different maps
of bundle and opportunity to define an operation over sections
create problems for this structure. 

			\section{Representation of Fibered \texorpdfstring{$\mathcal{F}$}{F}\Hyph Algebra}
			\label{section: Representation of Fibered F-Algebra}

	\begin{definition}
	\label{definition: nonsingular transformation, bundle} 
We call the fibered map
\[
t:\mathcal{E}\rightarrow\mathcal{E}
\]
\AddIndex{transformation of bundle}
{transformation of bundle},
if respective fiber map
\[
t_x:E_x\rightarrow E_x
\]
is transformation of a fiber.
	 \qed
	 \end{definition}

		\begin{theorem}
		\label{theorem: continuity of transformation}
Let $U$ be open set of base $M$ such
that there exists
a local chart of the bundle $\bundle{}pE{}$.
Let $t$ be transformation of bundle $\bundle{}pE{}$.
Let $x\in M$ and $p'=t_x(p)$ in fiber $E_x$.
Then for an open set $V\subseteq M$,
$x\in V$ and for an open set $W'\subseteq E$,
$p'\in W'$
there exists an open set $W\subseteq E$
such that if $x_1\in V$, $p_1\in W$,
then $p'_1=t_{x_1}(p_1)\in W$.
		\end{theorem}
		\begin{proof}
According to \citeBib{Bourbaki: General Topology 1}, page 44,
sets $V\times W$, where $V$ forms base
of topology of space $U$ and $W$ forms base
of topology of space $E$,
form base of topology of space $\mathcal{E}$. 

Since map $t$ is continuous, then for open set
$V\times W'$ there exists open set $V\times W$ such,
that $t( V\times W)\subseteq V\times W'$. This is the statement of theorem.
		\end{proof}

		\begin{theorem}
Transformation of bundle $\bundle{}pE{}$
maps section into section.
		\end{theorem}
		\begin{proof}
We define the image of section $s$ over transformation $t$
using commutative diagram
	\[
\xymatrix{
\mathcal{E}\ar[rr]^t&&\mathcal{E}\\
&M\ars[ul]^s\ars[ur]_{s'}&
}
	\]
Continuity of map $s'$
follows from theorem
\ref{theorem: continuity of transformation}.
		\end{proof}

	 \begin{definition}
	\label{definition: Tstar transformation of bundle} 
Transformation of bundle is
\AddIndex{left-side transformation}
{left-side transformation of bundle} or
\AddIndex{\Ts transformation of bundle}{Tstar transformation of bundle}
if it acts from left
\[
u'=t u
\]
We denote
\symb{{}^\star \mathcal{E}}1
{set of Tstar nonsingular transformations of bundle} or
\symb{{}^\star \bundle{}pE{}}1
{set of Tstar nonsingular transformations of bundle, projection} or
the set of nonsingular \Ts transformations of bundle
$\bundle{}pE{}$.
	 \qed
	 \end{definition}

	 \begin{definition}
Transformations is
\AddIndex{right-side transformations}
{right-side transformation of bundle} or
\AddIndex{\sT transformation of bundle}{starT transformation of bundle}
if it acts from right
\[
u'= ut
\]
We denote
\symb{\mathcal{E}^\star}1{set of starT nonsingular transformations of bundle} or
\symb{\bundle{}pA{}^\star}1{set of starT nonsingular transformations of bundle, projection}
the set of nonsingular \sT transformations of bundle
$\bundle{}pE{}$.
	 \qed
	 \end{definition}

We denote
\symb{e}1{identical transformation of bundle}
identical transformation of bundle.

Since we define \Ts transformation of bundle by fiber, then set ${}^\star \bundle{}pE{}$
is bundle isomorphic to the bundle $\bundle{}p{{}^\star E}{}$.

		\begin{definition}
		\label{definition: Tstar representation of fibered F-algebra} 
Suppose we defined the structure of fibered $\mathcal{F}$\Hyph algebra
on the set ${}^\star \bundle{}pA{}$
(\citeBib{Burris Sankappanavar}).
Let $\bundle{}pB{}$ be fibered $\mathcal{F}$\Hyph algebra.
We call homomorphism of fibered $\mathcal{F}$\Hyph algebras
	\begin{equation}
f:\bundle{}pB{}\rightarrow {}^\star \bundle{}pA{}
	\label{eq: Tstar representation of fibered F-algebra}
	\end{equation}
\AddIndex{left-side representation}
{left-side representation of fibered F-algebra} or
\AddIndex{\Ts representation of fibered $\mathcal{F}$\Hyph algebra $\bundle{}pB{}$}
{Tstar representation of fibered F-algebra}.
		\qed
		\end{definition}

		\begin{definition}
		\label{definition: starT representation of fibered F-algebra} 
Suppose we defined the structure of fibered $\mathcal{F}$\Hyph algebra
on the set $\bundle{}pA{}^\star$
(\citeBib{Burris Sankappanavar}).
Let $\bundle{}pB{}$ be fibered $\mathcal{F}$\Hyph algebra.
We call homomorphism of fibered $\mathcal{F}$\Hyph algebras
	\[
f:\bundle{}pB{}\rightarrow \bundle{}pA{}^\star
	\]
\AddIndex{right-side representation}
{right-side representation of fibered F-algebra} or
\AddIndex{\sT representation of fibered $\mathcal{F}$\Hyph algebra $\bundle{}pB{}$}
{starT representation of fibered F-algebra}.
		\qed
		\end{definition}
We extend to bundle representation theory convention
described in remark
\xRef{0701.238}{remark: left and right matrix notation}.
We can write duality principle in the following form

		\begin{theorem}[duality principle]
		\label{theorem: duality principle, fibered F-algebra representation}
Any statement which holds for \Ts representation of fibered $\mathcal{F}$\Hyph algebra $\bundle{}pA{}$
holds also for \sT representation of fibered $\mathcal{F}$\Hyph algebra $\bundle{}pA{}$.
		\end{theorem}

There are two ways to define a \Ts representation of $\mathcal{F}$\Hyph algebra $B$ in the bundle
$\bundle{}pA{}$. We can define or \Ts representation in the fiber,
either define \Ts representation in the set $\Gamma(\bundle{}pA{})$.
In the former case the representation defines the same transformation in
all fibers. In the later case the picture is less restrictive, however we do not have the whole
picture of the diversity of representations in the bundle. Studying the representation of
the fibered $\mathcal{F}$\Hyph algebra, we point out that representations in different fibers are
independent. Demand of smooth dependence of transformation on fiber
put additional constrains for \Ts representation of fibered $\mathcal{F}$\Hyph algebra.
The same time this constrain allows learn \Ts representation of the fibered $\mathcal{F}$\Hyph algebra
when in the fiber there defined $\mathcal{F}$\Hyph algebra with parameters
(for instance, the structure constants of Lie group) smooth dependent on fiber.

		\begin{remark}
		\label{remark: reduce details on the diagram, fibered F-algebra}
Using diagrams we can express definition
\ref{definition: Tstar representation of fibered F-algebra}
the following way.
	
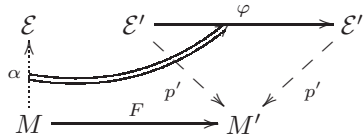
\[
\xymatrix{
\mathcal{E}&\mathcal{E}'\arp[dr]_{p'}\ar[rr]^\varphi="C2"
& & \mathcal{E}'\arp[dl]^{p'}\\
M\ars[u]^\alpha="C1"\ar[rr]^F& & M'\\
\ar @{=>}@/_1pc/ "C1";"C2"
}
	\]
Map $F$ is injection. Because we expect that representation of fibered $\mathcal{F}$\Hyph algebra
acts in each fiber, then we see that map $F$ is bijection.
Without loss of generality, we assume that $M=M'$ and map $F$ is the identity
map. We tell that we define the representation  of the fibered $\mathcal{F}$\Hyph algebra
$\bundle{}pB{}$ in the bundle $\bundle{}pA{}$ over the set $M$.
Since we know the base of the bundle, then to reduce details on the diagram we will
describe the representation using the following diagram
	\[
\xymatrix{
\mathcal{E}&\mathcal{E}'\ar[rr]^\varphi="C2"
& & \mathcal{E}'\\
M\ars[u]^\alpha="C1"\\
\ar @{=>}@/_1pc/ "C1";"C2"
}
	\]
		\qed
		\end{remark}

	 \begin{definition}
	 \label{definition: effective representation of fibered F-algebra}
Suppose map \eqref{eq: Tstar representation of fibered F-algebra} is
an isomorphism of the fibered $\mathcal{F}$\Hyph algebra $\bundle{}pB{}$ into ${}^\star \bundle{}pA{}$.
Then the \Ts representation of the fibered $\mathcal{F}$\Hyph algebra $\bundle{}pB{}$ is called
\AddIndex{effective}{effective representation of fibered F-algebra}.
	 \qed
	 \end{definition}

		\begin{remark}
		\label{remark: notation for effective representation of fibered F-algebra}
Suppose the \Ts representation of fibered $\mathcal{F}$\Hyph algebra is effective. Then we identify
an element of fibered $\mathcal{F}$\Hyph algebra and its image and write \Ts transformation
caused by element $a\in A$
as
\[v'=av\]
Suppose the \sT representation of $\mathcal{F}$\Hyph algebra is effective. Then we identify
an element of fibered $\mathcal{F}$\Hyph algebra and its image and write \sT transformation
caused by element $a\in A$
as
\[v'=va\]
		\qed
		\end{remark}

	 \begin{definition}
	 \label{definition: transitive representation of fibered F-algebra}
We call a \Ts representation of fibered $\mathcal{F}$\Hyph algebra
\AddIndex{transitive}{transitive representation of fibered F-algebra}
if for any $a, b \in V$ exists such $g$ that
\[a=f(g)b\]
We call a \Ts representation of fibered $\mathcal{F}$\Hyph algebra
\AddIndex{single transitive}{single transitive representation of fibered F-algebra}
if it is transitive and effective.
	 \qed
	 \end{definition}

		\begin{theorem}	
\Ts representation is single transitive if and only if for any $a, b \in M$
exists one and only one $g\in A$ such that $a=f(g)b$
		\end{theorem}
		\begin{proof}
Colorary of definitions \ref{definition: effective representation of fibered F-algebra}
and \ref{definition: transitive representation of fibered F-algebra}.
		\end{proof}

\def\texFiberedGroup{}

			\section{Representation of fibered group}
			\label{section: Representation of fibered group}

	\begin{definition}
	\label{definition: homomorphism of fibered groups} 
Let $\bundle{\mathcal{G}}{p}GM$
and $\bundle{\mathcal{G}'}{p'}{G'}{M'}$
be two fibered groups. Bundle map
\[
f:\mathcal{G}\rightarrow\mathcal{G}'
\]
is called
\AddIndex{homomorphism of fibered groups}{homomorphism of fibered groups}
if respective fiber map
\[
f_x:G_x\rightarrow G_{x'}'
\]
is homomorphism of groups.
	 \qed
	 \end{definition}

	\begin{definition}
	\label{definition: antihomomorphism of fibered groups} 
Let $\bundle{\mathcal{G}}{p}GM$
and $\bundle{\mathcal{G}'}{p'}{G'}{M'}$
be two fibered groups. Bundle map
\[
f:\mathcal{G}\rightarrow\mathcal{G}'
\]
is called
\AddIndex{antihomomorphism of fibered groups}{antihomomorphism of fibered groups}
if respective fiber map
\[
f_x:G_x\rightarrow G_{x'}'
\]
is antihomomorphism of groups.
	 \qed
	 \end{definition}

		\begin{definition}
		\label{definition: Tstar representation of fibered group}
Let $\bundle{}pG{}$ be fibered group.
We call map
	\begin{equation}
f:\bundle{}pG{}\rightarrow {}^\star \bundle{}pA{}
	\label{eq: Tstar representation of fibered group}
	\end{equation}
\AddIndex{\Ts representation of fibered group}
{Tstar representation of fibered group}
$\bundle{}pG{}$ in bundle $\bundle{}pA{}$ if map $f$ holds
	\begin{equation}
f(ab)\mu=f(a)(f(b)\mu)
	\label{eq: Tstar product of bundle transformations}
	\end{equation}
	\begin{equation}
f(\epsilon)=e
	\label{eq: Tstar identical transformation of bundle}
	\end{equation}
		 \qed
		 \end{definition}

		\begin{definition}
		\label{definition: starT representation of fibered group} 
Let $\bundle{}pG{}$ be fibered group.
We call map
	\begin{equation}
f:\bundle{}pG{}\rightarrow\bundle{}pA{}^\star
	\label{eq: right-side representation of fibered group}
	\end{equation}
\AddIndex{\sT representation of fibered group}
{starT representation of fibered group}
$\bundle{}pG{}$ in bundle $\bundle{}pA{}$ if map $f$ holds
	\begin{equation}
\mu f(ab)=(\mu f(a))f(b)
	\label{eq: starT product of bundle transformations}
	\end{equation}
	\begin{equation}
f(\epsilon)=e
	\label{eq: starT identical transformation of bundle}
	\end{equation}
		 \qed
		 \end{definition}

			\begin{theorem}
			\label{theorem: inverse bundle transformation}
For any $a\in \bundle{}pG{}$
	\begin{equation}
	\label{eq: inverse bundle transformation}
f(a^{-1})=f(a)^{-1}
	\end{equation}
		\end{theorem}
		\begin{proof}
Since \eqref{eq: Tstar product of bundle transformations}
and \eqref{eq: Tstar identical transformation of bundle}, we have
\[
\mu=e \mu=f(aa^{-1})\mu
=f(a)(f(a^{-1})\mu)
\]
This completes the proof.
		\end{proof}

		 \begin{theorem}
		 \label{theorem: Tstar covariant representation of fibered group}
Let ${}^\star \bundle{}pA{}$ be a fibered group with respect to multiplication
	\begin{equation}
(t_1t_2)\mu=t_1(t_2\mu)
	\label{eq: Tstar covariant product of bundle transformations}
	\end{equation}
and $e$ be unit of group ${}^\star \bundle{}pA{}$.
Let map \eqref{eq: Tstar representation of fibered group}
be a homomorphism of fibered group
	\begin{equation}
f(ab)=f(a)f(b)
	\label{eq: Tstar covariant representation of fibered group}
	\end{equation}
Then this map is representation of fibered group $\bundle{}pG{}$ which we call
\AddIndex{covariant \Ts representation of fibered group}
{Tstar covariant representation of fibered group}.
		\end{theorem}
		\begin{proof}
Since $f$ is homomorphism of fibered group,
we have $f(\epsilon)=e$.

Since
\eqref{eq: Tstar covariant product of bundle transformations}
and \eqref{eq: Tstar covariant representation of fibered group}, we have
\[
f(ab)\mu=(f(a)f(b))\mu
=f(a)(f(b)\mu)
\]

According definition \ref{definition: Tstar representation of fibered group}
$f$ is representation of fibered group.
		\end{proof}

We use following diagram to represent covariant \Ts representation of fibered group on the bundle
	\[
\xymatrix{
&&&&&\mathcal{E}\ar[dd]^{f(ab)}="f(ab)"\ar[lld]_{f(a)}="f(a)"\\
&&&\mathcal{E}\ar[drr]_{f(b)}="f(b)"&&\\
\mathcal{G}^2\ar[rr]&&\mathcal{G}&&&\mathcal{E}\\
\odot\ar@{=>}@/^1pc/[rr] &&&&&\\
M\ars@/^2pc/[uu]^a="a"\ar@/_2pc/[uu]_b="b"\ars[rr]^{id}&&
M\ars[uu]_{ab}="ab"&&&
\ar @{=>}@/_4pc/ "ab";"f(ab)"
\ar @{=>}@/^4pc/ "a";"f(a)"
\ar @{=>}@/^3pc/ "b";"f(b)"
}
	\]

	 \begin{theorem}
	 \label{theorem: Tstar contravariant representation of fibered group}
Let ${}^\star \bundle{}pA{}$ be a fibered group with respect to multiplication
	\begin{equation}
(t_2t_1)\mu=t_1(t_2\mu)
	\label{eq: Tstar contravariant product of bundle transformations}
	\end{equation}
and $e$ be unit of fibered group ${}^\star \bundle{}pA{}$.
Let map \eqref{eq: Tstar representation of fibered group}
be an antihomomorphism of fibered group
	\begin{equation}
f(ba)=f(a)f(b)
	\label{eq: Tstar contravariant representation of fibered group}
	\end{equation}
Then this map is representation of fibered group $\bundle{}pG{}$ which we call
\AddIndex{contravariant \Ts representation of fibered group}
{Tstar contravariant representation of fibered group}.
		\end{theorem}
		\begin{proof}
Since $f$ is antihomomorphism of fibered group,
we have $f(\epsilon)=e$.

Since
\eqref{eq: Tstar contravariant product of bundle transformations}
and \eqref{eq: Tstar contravariant representation of fibered group}, we have
\[
f(ab)\mu=(f(b)f(a))\mu
=f(a)(f(b)\mu)
\]

According definition \ref{definition: Tstar representation of fibered group}
$f$ is representation of fibered group.
		\end{proof}

	\begin{example}
The group composition on fibered group
determines two different presentations on the fibered group:
the \AddIndex{\Ts shift on the fibered group}{Tstar shift, fibered group}
which we introduce by the equation
\symb{a\star b}0{Tstar shift, fibered group}
	\begin{equation}
b'=\ShowSymbol{Tstar shift, fibered group}=ab
	\label{eq: Tstar shift, fibered group}
	\end{equation}
and the \AddIndex{\sT shift on fibered group}{starT shift, fibered group}
which we introduce by the equation
\symb{b\star a}0{starT shift, fibered group}
	\begin{equation}
b'=\ShowSymbol{starT shift, fibered group}=ba
	\label{eq: starT shift, fibered group}
	\end{equation}
	 \qed
	 \end{example}

	\begin{example}
Let $\bundle{}p{GL}{}$ be bundle over set of real numbers.
Given the matrix $A$, we can define section $a(t)=exp(tA)$,
and this section will cause respective \Ts shift.
	 \end{example}

	 \begin{definition}
Let $f$ be representation of fibered group
$\bundle{}pG{}$ in bundle $\bundle{}pA{}$.
For any cut $u$ of bundle $\bundle{}pA{}$ we define its
\AddIndex{orbit of representation of fibered group}
{orbit of representation of fibered group} as set
\symb{\mathcal{O}(u,g\in\Gamma(\bundle{}pG{}),f(g)u)}
0{orbit of representation of fibered group}
\[
\ShowSymbol{orbit of representation of fibered group}
=\{v=f(g)u:g\in\Gamma(\bundle{}pG{})\}
\]
	 \qed
	 \end{definition}

Since $f(\epsilon)=e$ we have
$u\in\mathcal{O}(u,g\in\Gamma(\bundle{}pG{}),f(g)u)$.

		\begin{theorem}
		\label{theorem: proper definition of orbit, fibered group}
Suppose
	\begin{equation}
	\label{eq: orbit, fibered group, proposition}
v\in\mathcal{O}(u,g\in\Gamma(\bundle{}pG{}),f(g)u)
	\end{equation}
Then
\[
\mathcal{O}(u,g\in\Gamma(\bundle{}pG{}),f(g)u)
=\mathcal{O}(v,g\in\Gamma(\bundle{}pG{}),f(g)v)
\]
		\end{theorem}
		\begin{proof}
From \eqref{eq: orbit, fibered group, proposition} it follows
that there exists $\mu\in\Gamma(\bundle{}pG{})$ such that
	\begin{equation}
v=f(\mu)u
	\label{eq: orbit, fibered group, 1}
	\end{equation}
Suppose $\delta\in\mathcal{O}(v,g\in\Gamma(\bundle{}pG{}),f(g)v)$. Then
there exists $\nu\in\Gamma(\bundle{}pG{})$ such that
	\begin{equation}
\delta=f(\nu)v
	\label{eq: orbit, fibered group, 2}
	\end{equation}
If we substitude \eqref{eq: orbit, fibered group, 1}
into \eqref{eq: orbit, fibered group, 2} we get
	\begin{equation}
\delta=f(\nu)f(\mu)u
	\label{eq: orbit, fibered group, 3}
	\end{equation}
Since \eqref{eq: Tstar product of bundle transformations}
we see that
from \eqref{eq: orbit, fibered group, 3} it follows
that $\delta\in\mathcal{O}(u,g\in\Gamma(\bundle{}pG{}),f(g)u)$.
Thus
\[
\mathcal{O}(v,g\in\Gamma(\bundle{}pG{}),f(g)v)
\subseteq\mathcal{O}(u,g\in\Gamma(\bundle{}pG{}),f(g)u)
\]

Since \eqref{eq: inverse bundle transformation},
we see that
from \eqref{eq: orbit, fibered group, 1} it follows
that
	\begin{equation}
u=f(\mu)^{-1}v=f(\mu^{-1})v
	\label{eq: orbit, fibered group, 4}
	\end{equation}
From \eqref{eq: orbit, fibered group, 4} it follows that
$u\in\mathcal{O}(v,g\in\Gamma(\bundle{}pG{}),f(g)v)$ and therefore
\[
\mathcal{O}(u,g\in\Gamma(\bundle{}pG{}),f(g)u)
\subseteq\mathcal{O}(v,g\in\Gamma(\bundle{}pG{}),f(g)v)
\]
This completes the proof.
		\end{proof}

Let us define the representation of group $G$ on the bundle $$\bundle{\mathcal{E}}pA{M}$$
Since we call the representation transitive, then orbit of a point is
the manifold $\mathcal{E}$. In the case of
representation of fibered group $\bundle{}pG{}$ the orbit of a point is the fiber the point belongs to.

		 \begin{theorem}
		 \label{theorem: direct product of representations of fibered group}
Suppose $f_1$ is representation of fibered group $\bundle{}pG{}$
in bundle $\bundle{}p{A_1}{}$ and
$f_2$ is representation of fibered group $\bundle{}pG{}$
in bundle $\bundle{}p{A_2}{}$.
Then we introduce
\AddIndex{direct product of representations
$f_1$ and $f_2$ of fibered group}
{direct product of representations of fibered group} $\bundle{}pG{}$
	\begin{align*}
f&=f_1\otimes f_2:
\bundle{}pG{}\rightarrow \bundle{}p{A_1}{}\otimes \bundle{}p{A_2}{}\\
f(g)&=(f_1(g),f_2(g))
	\end{align*}
		\end{theorem}
		\begin{proof}
To show that $f$ is representation
it is enough to prove that $f$ satisfy to definition
\ref{definition: Tstar representation of fibered group}.
\[f(e)=(f_1(e),f_2(e))
=(e_1,e_2)=e\]
	\begin{align*}
f(ab)u&=(f_1(ab)u_1,f_2(ab)u_2)\\
&=(f_1(a)(f_1(b)u_1),
f_2(a)(f_2(b)u_2))\\
&=f(a)(f_1(b)u_1,f_2(b)u_2)\\
&=f(a)(f(b)u)
	\end{align*}
		\end{proof}

			\section{Single Transitive Representation}

	 \begin{definition}
	 \label{definition: kernel of inefficiency of representation of fibered group}
We call
\AddIndex{kernel of inefficiency of representation of fibered group}
{kernel of inefficiency of representation of fibered group} $\bundle{}pG{}$
a set \[K_f=\{g\in\Gamma(\bundle{}pG{}):f(g)=e\}\]
If $K_f=\{e\}$ we call representation of fibered group $G$
\AddIndex{effective}{effective representation of fibered group}.
	 \qed
	 \end{definition}

		\begin{theorem}	
		\label{theorem: kernel of inefficiency, representation of fibered group}
A kernel of inefficiency is a subgroup of fibered group $G$.
		\end{theorem}
		\begin{proof}
The proof does not depend on whether we use covariant representation or
contravariant representation. Assume $f$ is covariant representation and
$f(a_1)=\delta$ and $f(a_2)=\delta$. Then
\[f(a_1a_2)=f(a_1)f(a_2)=\delta\]
\[f(a^{-1})=f^{-1}(a)=\delta\]
		\end{proof}


		\begin{theorem}	
Representation is single transitive if and only if for any $a, b \in\Gamma(\bundle{}pA{})$
exists one and only one $g\in\bundle{}pG{}$ such that $a=f(g)b$
		\end{theorem}
		\begin{proof}
Statement is colorary of definitions
\ref{definition: effective representation of fibered F-algebra} and
\ref{definition: kernel of inefficiency of representation of fibered group} and
of the theorem \ref{theorem: kernel of inefficiency, representation of fibered group}.
		\end{proof}

	 \begin{definition}
We call a bundle $\bundle{}pA{}$
\AddIndex{homogeneous bundle of fibered group}
{homogeneous bundle of fibered group} $\bundle{}pG{}$
if we have single transitive representation of fibered group $\bundle{}pG{}$ on $\bundle{}pA{}$.
	 \qed
	 \end{definition}

		\begin{theorem}	%
		\label{theorem: single transitive representation of fibered group}
If we define a single transitive \Ts representation $f$ of the fibered group $\bundle{}pG{}$
on the bundle $\bundle{}pA{}$
then we can uniquely define coordinates on $\bundle{}pA{}$
using coordinates on the fibered group $\bundle{}pG{}$.

If $f$ is a covariant single transitive representation
than $f(a)$ is equivalent to the left shift $a\star$ on the fibered group $\bundle{}pG{}$.
If $f$ is a contravariant single transitive representation than
$f(a)$ is equivalent to the right shift $\star a$ on the fibered group $\bundle{}pG{}$.
		\end{theorem}
		\begin{proof}
The representation of the fibered group $\bundle{}pG{}$ in the bundle
$\bundle{}pA{}$ is single transitive
iff the representation of the group $G$ in the fiber $A_x$ for any $x$ is single transitive.
Let representation $f$ of the fibered group $\bundle{}pG{}$ be 
a covariant single transitive representation.
Let $u$, $v$ be sections of the bundle $\bundle{}pA{}$.
and $x\in M$. According to theorem
\xRef{0701.238}{theorem: single transitive representation of group}
we get the only element $g_x\in G$ such that
\[
u(x)=f(g_x)v(x)=g_xv(x)
\]a contravariant single transitive representation
Thus the map $x\rightarrow g_x$ is the section of the fibered group $\bundle{}pG{}$.

The same way we prove the statement for a covariant single transitive representation.

To prove the first statement, we need to select
the map on the manifold $M$, where both bundles are trivial. Then we can
represent coordinates of the point $u\in\bundle{}pA{}$ as tuple of coordinates
$(x,y)$ where $x$ are coordinates of projection to the manifold $M$ and $y$ are coordinates
of the point in the fiber. We can
represent coordinates of the point $a\in\bundle{}pG{}$ as tuple of coordinates
$(x,g)$ where $x$ are coordinates of projection to the manifold $M$ and $g$ are coordinates
of the point in the group. Respectively, coordinates of the section of the bundle $\bundle{}pA{}$ are
the map $y:M\rightarrow A$, and coordinates of the section of the bundle $\bundle{}pG{}$ are
the map $y:M\rightarrow G$.

We select a section $v\in\Gamma(\bundle{}pA{})$
and define coordinates of a point $w\in\Gamma(\bundle{}pA{})$
as coordinates of $a\in\bundle{}pG{}$ such that $w=f(a) v$.
Coordinates defined this way are unique
up to choice of an initial section $v\in\Gamma(\bundle{}pA{})$
because the action is effective.
		\end{proof}

		\begin{remark}
We will write effective \Ts covariant representation of the fibered group $\bundle{}pG{}$ as
\[v'=a\star v=av\]
Orbit of this representation is
\symb{\bundle{}pG{}v}0{orbit of effective Tstar covariant representation of fibered group}
\[
\ShowSymbol{orbit of effective Tstar covariant representation of fibered group}
=\bundle{}pG{}\star v
\]
		\qed
		\end{remark}

		\begin{remark}
We will write effective \sT covariant representation of the fibered group $\bundle{}pG{}$ as
\[v'=v\star a=va\]
Orbit of this representation is
\symb{v\bundle{}pG{}}0{orbit of effective starT covariant representation of fibered group}
\[
\ShowSymbol{orbit of effective starT covariant representation of fibered group}
=v\star\bundle{}pG{}
\]
		\qed
		\end{remark}

		\begin{theorem}	%
		\label{theorem: shifts on fibered group commuting}
Left and right shifts on the fibered group $\bundle{}pG{}$ are commuting.
		\end{theorem}
		\begin{proof}
This is the consequence of the associativity on the fibered group $\bundle{}pG{}$
\[(a\star \star b)c = a(cb)=(ac)b=(\star b\ a\star)c\]
		\end{proof}

		\begin{theorem}	%
		\label{theorem: two representations of fibered group}
If we defined a single transitive covariant \Ts representation $f$
of the fibered group $\bundle{}pG{}$ on the bundle $\bundle{}pA{}$
then we can uniquely define a single transitive covariant \sT representation $h$
of the fibered group $\bundle{}pG{}$ on the bundle $\bundle{}pA{}$
such that diagram
	
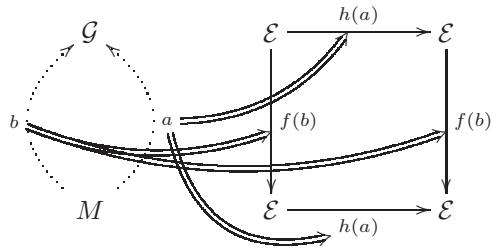
\[
\xymatrix{
\mathcal{G}&&\mathcal{E}\ar[rr]^{h(a)}="Ct2"\ar[dd]^{f(b)}="Cl2"
& & \mathcal{E}\ar[dd]^{f(b)}="Cr2"\\
&&&&\\
M\ars@/^2pc/[uu]^b="Cl1"\ars@/_2pc/[uu]_a="Cr1"&&\mathcal{E}\ar[rr]_{h(a)}="Cb2"
& & \mathcal{E}\\
\ar @{=>}@/_1pc/ "Cr1";"Ct2"
\ar @{=>}@/_2pc/ "Cr1";"Cb2"
\ar @{=>}@/_1pc/ "Cl1";"Cl2"
\ar @{=>}@/_1.5pc/ "Cl1";"Cr2"
}
	\]
is commutative for any $a$, $b\in \Gamma(\bundle{}pG{})$.
		\end{theorem}
		\begin{proof}
Let $f$ be a single transitive covariant \Ts representation.
In each fiber $A_x$ the representation $f$ defines a single transitive covariant \Ts representation
$f_x$ of group $G$. According to theorem
\xRef{0701.238}{theorem: two representations of group}
in fiber $A_x$ we uniquely define a single transitive covariant \sT representation
$h_x$ comutable with representation $f_x$.
For a section $a\in\Gamma(\bundle{}pG{})$ we define the section
\[
h(a):x\rightarrow h_x(a)
\] 
Map $h$ is homomorphism of fibered group.
		\end{proof}

We call representations $f$ and $h$
\AddIndex{twin representations of the fibered group}{twin representations of fibered group}
$\bundle{}pG{}$.
\OpenBiblio


\BiblioItem{texIntro}{Einstein: Geometry and Experience}
{
Einstein, Geometry and Experience, (1921)
}%

\BiblioItem{texGenRelativity}{Ghez}
{
Ghez et al.,
The First Measurement of Spectral Lines in a Short-Period Star Bound to the Galaxy's Central Black Hole: A Paradox of Youth,
\href{http://www.journals.uchicago.edu/ApJ/journal/issues/ApJL/v586n2/16990/brief/16990.abstract.html}{ApJL, 586, L127} (2003),
eprint \href{http://arxiv.org/abs/astro-ph/0302299}{arXiv:astro-ph/0302299} (2003)
}%

\BiblioItem{texGenRelativity}{Schodel}
{
R. Sch\"odel et al.,
A star in a 15.2-year orbit around the supermassive black hole at the centre of the Milky Way,
\href{http://www.nature.com/cgi-taf/DynaPage.taf?file=/nature/journal/v419/n6908/abs/nature01121_fs.html}{Nature 419, 694} (2002)
}%

\BiblioItem{texAffine}{Mielke}
{
Eckehard W. Mielke, Affine generalization of the Komar complex of general relativity,
\href{http://prola.aps.org/searchabstract/PRD/v63/i4/e044018}{Phys. Rev. D 63, 044018} (2001)
}%

\BiblioItem{texAffine}{Obukhov}
{
Yu. N. Obukhov and J. G. Pereira, Metric-affine approach to teleparallel gravity,
\href{http://scitation.aip.org/getabs/servlet/GetabsServlet?prog=normal&id=PRVDAQ000067000004044016000001&idtype=cvips&gifs=Yes}
{Phys. Rev. D 67, 044016} (2003),
eprint \href{http://arxiv.org/abs/gr-qc/0212080}{arXiv:gr-qc/0212080} (2002)
}%

\BiblioItem{texAffine}{Sardanashvily}
{
Giovanni Giachetta, Gennadi Sardanashvily, Dirac Equation in Gauge and Affine-Metric Gravitation Theories,
eprint \href{http://arxiv.org/abs/gr-qc/9511035}{arXiv:gr-qc/9511035} (1995)
}%

\BiblioItem{texAffine}{Gauge}
{
Frank Gronwald and Friedrich W. Hehl, On the Gauge Aspects of Gravity, eprint
\href{http://arxiv.org/abs/gr-qc/9602013}{arXiv:gr-qc/9602013} (1996)
}%

\BiblioItem{texAffine}{Neeman}
{
Yuval Neeman, Friedrich W. Hehl, Test Matter in a Spacetime with Nonmetricity, eprint
\href{http://arxiv.org/abs/gr-qc/9604047}{arXiv:gr-qc/9604047} (1996)
}%

\BiblioItem{texAffine}{0405.027}
{
Aleks Kleyn,
Reference Frame in General Relativity,
eprint \href{http://arxiv.org/abs/gr-qc/0405027}{arXiv:gr-qc/0405027} (2004)
}%

\BiblioItem{texTidal,texAffine}{torsion}
{
F. W. Hehl, P. von der Heyde, G. D. Kerlick, and J. M. Nester,
General relativity with spin and torsion: Foundations and prospects,
\href{http://prola.aps.org/abstract/RMP/v48/i3/p393_1}{Rev. Mod. Phys. 48, 393} (1976)
}%

\BiblioItem{texTidal,texNewton}{Megged}
{
O. Megged, Post-Riemannian Merger of Yang-Mills Interactions with Gravity,
eprint \href{http://arxiv.org/abs/hep-th/0008135}{arXiv:hep-th/0008135} (2001)
}%


\BiblioItem{texNewton}{gr-qc-9604027}
{
Yu.N. Obukhov, E.J. Vlachynsky, W. Esser, R. Tresguerres and F.W. Hehl,
An exact solution of the metric-affine gauge theory with dilation, shear, and spin charges,
eprint \href{http://arxiv.org/abs/gr-qc/9604027}{arXiv:gr-qc/9604027} (1996)
}%

\BiblioItem{texLagrange}{Weinberg}
{
Steven Weinberg. The Quantum Theory of Fields. Cambridge university press.
}%

\BiblioItem{texLagrange}{Reinhardt}
{
Greiner Reinhardt. Field Quantization. Springer.
}%

\BiblioItem{texLagrange}{Landau}
{
L. D. Landau, E. M. Lifshich, The classical theory of fields.
Oxford, New York, Pergamon Press
}%

\BiblioItem{texTidal}{Wheeler}
{
Ignazio Ciufolini, John Wheeler. Gravitation and Inertia.
Princeton university press.
}%

\BiblioItem{texTidal}{0405.028}
{
Aleks Kleyn, Metric-Affine Manifold,
eprint \href{http://arxiv.org/abs/gr-qc/0405028}{arXiv:gr-qc/0405028} (2004)
}%

\BiblioItem{texTidal}{Anderson02}
{
J. D. Anderson, P. A. Laing, E. L. Lau, A. S. Liu, M. M. Nieto, and S. G. Turyshev,
Study of the anomalous acceleration of Pioneer 10 and 11,
\href{http://prola.aps.org/searchabstract/PRD/v65/i8/e082004}{Phys. Rev. D 65, 082004, 50 pp.}, (2002),
eprint \href{http://arxiv.org/abs/gr-qc/0104064}{arXiv:gr-qc/0104064} (2001)
}%

\BiblioItem{texTidal}{Anderson98}
{
J. D. Anderson, P. A. Laing, E. L. Lau, A. S. Liu, M. M. Nieto, and S. G. Turyshev,
Indication, from Pioneer 10/11, Galileo, and Ulysses Data, of an Apparent Anomalous, Weak, Long-Range Acceleration,
\href{http://prola.aps.org/abstract/PRL/v81/i14/p2858_1}{Phys. Rev. Lett. 81, 2858}, (1998),
eprint \href{http://arxiv.org/abs/gr-qc/9808081}{arXiv:gr-qc/9808081} (1998)
}%


\BiblioItem{texReferenceFrame,texFiberedAlgebra}{Serge Lang}
{
Serge Lang,
Algebra, Springer, 2002
}%

\BiblioItem{texFiberedAlgebra,texTstarMorphism}{Burris Sankappanavar}
{
S. Burris, H.P. Sankappanavar,
A Course in Universal Algebra, Springer-Verlag (March, 1982),
\\eprint
\href{http://www.math.uwaterloo.ca/~snburris/htdocs/ualg.html}
{http://www.math.uwaterloo.ca/~snburris/htdocs/ualg.html}
\\(The Millennium Edition)
}%


\BiblioItem{texAffine,texRepresentation,texBasis,texDrcBasis,texLinearMap,texVectorSpace}{Rashevsky}
{
P. K. Rashevsky, Riemann Geometry and Tensor Calculus,\\
Moscow, Nauka, 1967
}%

\BiblioItem{texDrcBasis,texBasis}{Korn}
{
Granino A. Korn, Theresa M. Korn,
Mathematical Handbook for Scientists and Engineer,
McGraw-Hill Book Company, New York, San Francisco,
Toronto, London, Sydney, 1968
}%


\BiblioItem{texGenRelativity}{Tartaglia}
{
Angelo Tartaglia and Matteo Luca Ruggiero,
Angular Momentum Effects in Michelson\Hyph Morley Type Experiments,
Gen.Rel.Grav. 34, 1371-1382 (2002),\\
eprint \href{http://arxiv.org/abs/gr-qc/0110015}{arXiv:gr-qc/0110015} (2001)
}%

\BiblioItem{texGenRelativity}{Tomozawa}
{
Yukio Tomozawa, Speed of Light in Gravitational Fields, eprint
\href{http://arxiv.org/abs/astro-ph/0303047}{arXiv:astro-ph/0303047} (2004)
}%

\BiblioItem{texGenRelativity}{Magueijo}
{
Joao Magueijo,
Covariant and locally Lorentz-invariant varying speed of light theories,
\href{http://prola.aps.org/abstract/PRD/v62/i10/e103521}{Phys. Rev. D 62, 103521} (2000),
eprint \href{http://arxiv.org/abs/gr-qc/0007036}{arXiv:gr-qc/0007036} (2000)
}%

\BiblioItem{texGenRelativity}{Bassett}
{
Bruce A. Bassett, Stefano Liberati, Carmen Molina-Paris, and Matt Visser,
Geometrodynamics of variable-speed-of-light cosmologies,
\href{http://prola.aps.org/abstract/PRD/v62/i10/e103518}{Phys. Rev. D 62}, 103518 (2000),
eprint \href{http://arxiv.org/abs/astro-ph/0001441}{arXiv:astro-ph/0001441} (2000)
}%

\BiblioItem{texGenRelativity}{Straumann}
{
Lochlainn O'Raifeartaigh and Norbert Straumann,
Gauge theory: Historical origins and some modern developments,
\href{http://prola.aps.org/abstract/RMP/v72/i1/p1_1}{Rev. Mod. Phys. 72, 1} (2000)
}%

\BiblioItem{texGenRelativity}{Lammerzahl}
{
Claus L\"ammerzahl, Mark P. Haugan,
On the interpretation of Michelson\Hyph Morley experiments,
{Phys. Lett. A282 223-229} (2001),\\
eprint \href{http://arxiv.org/abs/gr-qc/0103052}{arXiv:gr-qc/0103052} (2001)
}%

\BiblioItem{texGenRelativity}{Muller}
{
Holger Muller et al.,
Modern Michelson-Morley Experiment using Cryogenic Optical Resonators,
\href{http://prola.aps.org/searchabstract/PRL/v91/i2/e020401}{Phys. Rev. Lett. 91, 020401} (2003),
eprint \href{http://arxiv.org/abs/physics/0305117}{arXiv:physics/0305117} (2000)
}%

\BiblioItem{texGenRelativity,texTidal}{Ranada}
{
Antonio F. Ranada,
Pioneer acceleration and variation of light speed: experimental situation,
eprint \href{http://arxiv.org/abs/gr-qc/0402120}{arXiv:gr-qc/0402120} (2004)
}%

\BiblioItem{texBiring,texVectorSpace}{math.QA-0208146}
{
I. Gelfand, S. Gelfand, V. Retakh, R. Wilson,
Quasideterminants,\\
eprint \href{http://arxiv.org/abs/math.QA/0208146}{arXiv:math.QA/0208146} (2002)
}%

\BiblioItem{texBiring,texVectorSpace}{q-alg-9705026}
{
I.Gelfand, V.Retakh,
Quasideterminants, I,\\
eprint \href{http://arxiv.org/abs/q-alg/9705026}{arXiv:q-alg/9705026} (1997)
}%

\BiblioItem{texVectorSpace}{Gelfand Retakh 1991}
{
I. Gelfand and V. Retakh, Determinants of Matrices over Noncommutative Rings, Funct.
Anal. Appl. 25 (1991), no. 2, 91-102
}%

\BiblioItem{texVectorSpace}{Gelfand Retakh 1992}
{
I. Gelfand and V. Retakh, A Theory of Noncommutative Determinants and Characteristic
Functions of Graphs, Funct. Anal. Appl. 26 (1992), no. 4, 1-20
}%

\BiblioItem{texVectorSpace}{hep-th-9407124}
{
I. M. Gelfand, D. Krob, A. Lascoux, B. Leclerc, V.S. Retakh and J.-Y. Thibon,
Noncommutative symmetric functions,\\
eprint \href{http://arxiv.org/abs/hep-th/9407124}{arXiv:hep-th/9407124} (1994)
}%

\BiblioItem{texVectorSpace}{Carl Faith 1}
{
Carl Faith, Algebra: Rings, Modules and Categories I,
Springer - Verlag, Berlin - Heidelberg - New York, 1973
}%



\BiblioItem{texReferenceFrame}{math.DG-0412391}
{
Aleks Kleyn,
Basis Manifold,
eprint \href{http://arxiv.org/abs/math.DG/0412391}{arXiv:math.DG/0412391} (2004)
}%

\BiblioItem{texFiberedAlgebra,texBundleRelation,texTstarMorphism}{0701.238}
{
Aleks Kleyn,
Lectures on Linear Algebra over Skew Field,\\
eprint \href{http://arxiv.org/abs/math.GM/0701238}{arXiv:math.GM/0701238} (2007)
}%

\BiblioItem{texBundleRelation,texPrefaceRelation}{0702.561}
{
Aleks Kleyn,
Algebra Bundle,\\
eprint \href{http://arxiv.org/abs/math.DG/0702561}{arXiv:math.DG/0702561} (2007)
}%


\BiblioItem{texPolymodule}{math.RA-0501237v1}
{
Aleks Kleyn,
Module Over Skew-Field, version 1,\\
eprint \href{http://arxiv.org/abs/math/0501237v1}{arXiv:math.RA/0501237v1} (2005)
}%

\ifx\PrintBook\Defined
\else
\BiblioItem{texVectorSpace,texFiberedAlgebra}{0612.111}
{
Aleks Kleyn,
Biring of Matrices,\\
eprint \href{http://arxiv.org/abs/math.OA/0612111}{arXiv:math.OA/0612111} (2006)
}%
\fi

\BiblioItem{texHomotopy}{q-alg-9705009}
{
John C. Baez,
An Introduction to n-Categories,\\
eprint \href{http://arxiv.org/abs/q-alg/9705009}{arXiv:q-alg/9705009} (1997)
}%

\BiblioItem{texIntro}{Einstein: Isaak Newton}
{
Albert Einstein,
Isaak Newton,
Manchester Guardian, 16, 234 - 235 (1927)
}%

\BiblioItem{texPrefaceRelation}{Tolstoi about Anna Karenina}
{
Tolstoi about Anna Karenina,
in book A Karenina Companion, by C. J. G. Turner,
published by Wilfrid Laurier University Press (August 1993)
}%

\BiblioItem{texBundleRelation,texPrefaceRelation,texTstarMorphism,texCartesian}
{Cohn: Universal Algebra}
{
Paul M. Cohn,
Universal Algebra,
Springer, 1981
}%

\BiblioItem{texCartesian}
{Maunder: Algebraic Topology}
{
C. R. F. Maunder,
Algebraic Topology,
Dover Publications, Inc, Mineola, New York, 1996
}%

\BiblioItem{texFiberedAlgebra}{Pommaret: Partial Differential Equations}
{
J.-F. Pommaret,
Partial Differential Equations and Group Theory,
Springer, 1994
}%

\BiblioItem{texBundleRelation}{Bourbaki: Set Theory}
{
N. Bourbaki,
Theory of sets,
Springer, 2004
}%

\BiblioItem{texCalculus}{Pontryagin: Topological Group}
{
L. S. Pontryagin,
Selected Works, Volume Two, Topological Groups,
Gordon and Breach Science Publishers, 1986
}

\BiblioItem{texCartesian,texFiberedAlgebra}{Bourbaki: General Topology 1}
{
N. Bourbaki,
General Topology, Chapters 1 - 4,
Springer, 1989
}

\BiblioItem{texCalculus}{Bourbaki: Topological Vector Space}
{
N. Bourbaki,
Topological Vector Spaces, Chapters 1 - 5,
Springer, 1987
}

\BiblioItem{texCalculus}{Bourbaki: General Topology: Chapter 5 - 10}
{
N. Bourbaki,
General Topology, Chapters 5 - 10,
Springer, 1989
}

\BiblioItem{texFiberedAlgebra}{Hatcher: Algebraic Topology}
{
Allen Hatcher,
Algebraic Topology,
Cambridge University Press, 2002
}

\CloseBiblio

\OpenIndex
\SetIndexSpace%
\Index{texLinearMap}
   {$1$-\drc form}%
   {1-drc form, vector spaces}%
\SetIndexSpace%
\Index{texPolymodule}
   {$(2)$\hyph vector space}%
   {(2)-vector space}%
\Index{texBundleRelation}
   {$2$\Hyph ary fibered relation}%
   {2 ary fibered relation}%
\SetIndexSpace%
\Index{texCalculus}
   {$A$\Hyph valued function}%
   {A valued function}%
\Index{texBiring}
   {$(^{\gi a}_{\gi b})$\hyph \CR quasideterminant}%
   {a b cr-quasideterminant}%
\Index{texBiring}
   {$(^{\gi a}_{\gi b})$\hyph \RC quasideterminant}%
   {a b RC-quasideterminant}%
\Index{texCalculus}
   {absolute value on skew field}%
   {absolute value on skew field}%
\Index{texBasis}
   {active representation}%
   {active representation}%
\Index{texDrcBasis}
   {active \sT representation}%
   {active representation, vector space}%
\Index{texBasis}
   {active transformation on basis manifold}%
   {active transformation}%
\Index{texBasis}
   {affine basis}%
   {Affine Basis}%
\Index{texBasis}
   {affine transformation group}%
   {AffineTransformationGroup}%
\Index{texTypeBasis}
   {affine transformation group}%
   {AffineTransformationGroup}%
\Index{texBasis}
   {affine transformation on basis manifold}%
   {affine transformation}%
\Index{texBiring}
   {alternative representation of matrix}%
   {Alternative representation}%
\Index{texReferenceFrame}
   {anholonomic coordinate}%
   {anholonomic coordinate}%
\Index{texReferenceFrame}
   {anholonomic coordinates of vector}%
   {vector anholonomic coordinates}%
\Index{texReferenceFrame}
   {anholonomic coordinates on manifold}%
   {anholonomic coordinates on manifold}%
\Index{texReferenceFrame}
   {anholonomity object}%
   {anholonomity object}%
\Index{texFiberedGroup}
   {antihomomorphism of fibered groups}%
   {antihomomorphism of fibered groups}%
\Index{texBundleRelation}
   {antisymmetric $2$\Hyph ary fibered relation}%
   {antisymmetric 2 ary fibered relation}%
\Index{texFiberedAlgebra}
   {arity of operation}%
   {arity of operation}%
\Index{texTstarRepresentation}
   {associative law for covariant \Ts representation}%
   {associative law for Tstar covariant representation}%
\Index{texDrcMorphism}
   {associative law for \drc linear maps of vector spaces}%
   {associative law for drc linear maps of vector spaces}%
\Index{texVectorSpace}
   {associative law for \Ts vector space}%
   {associative law, Tstar vector space}%
\Index{texLinearMap}
   {associative law for twin representations}%
   {associative law for twin representations}%
\Index{texBundleRelation}
   {associative law of composition of fibered correspondences}%
   {associative law, composition of fibered correspondences}%
\Index{texAffine}
   {auto parallel line}%
   {auto parallel line}%
\SetIndexSpace%
\Index{texBundleRelation}
   {base of fibered correspondence}%
   {base of fibered correspondence}%
\Index{texCartesian}
   {base of map}%
   {base of map}%
\Index{texTypeBasis}
   {basis}%
   {}%
\SubIndex{texTypeBasis}
   {affine}%
   {Affine Basis}%
\SubIndex{texTypeBasis}
   {central affine}%
   {Central Affine Basis}%
\SubIndex{texTypeBasis}
   {orthonornal}%
   {Orthonornal Basis}%
\Index{texDrcBasis}
   {basis manifold}%
   {}%
\SubIndex{texTypeBasis}
   {of affine space}%
   {Basis Manifold, Affine Space}%
\SubIndex{texTypeBasis}
   {of central affine space}%
   {Basis Manifold, Central Affine Space}%
\SubIndex{texDrcBasis}
   {of \drc vector space}%
   {basis manifold of vector space}%
\SubIndex{texTypeBasis}
   {of Euclid space}%
   {Basis Manifold, Euclid Space}%
\Index{texBasis}
   {basis manifold of affine space}%
   {Basis Manifold, Affine Space}%
\Index{texBasis}
   {basis manifold of central affine space}%
   {Basis Manifold, Central Affine Space}%
\Index{texBasis}
   {basis manifold of Euclid space}%
   {Basis Manifold, Euclid Space}%
\Index{texBasis}
   {basis manifold of vector space}%
   {basis manifold of vector space}%
\Index{texBasis}
   {basis of vector space}%
   {Basis}%
\Index{texLieRepresentation}
   {basis vector}%
   {}%
\SubIndex{texLieRepresentation}
   {of \sT representation}%
   {basis vector of starT representation}%
\SubIndex{texLieRepresentation}
   {of \Ts representation}%
   {basis vector of Tstar representation}%
\Index{texBiring}
   {biring}%
   {biring}%
\SetIndexSpace%
\Index{texVectorSpace}
   {\subs rows \dcr vector space}%
   {subs rows dcr vector space}%
\Index{texBiring}
   {\subs row of matrix}%
   {c row}%
\Index{texBiring}
   {$c$\hyph row of matrix}%
   {c-row}%
\Index{texAffine}
   {Cartan connection}%
   {Cartan connection}%
\Index{texAffine}
   {Cartan curvature}%
   {Cartan curvature}%
\Index{texAffine}
   {Cartan derivative}%
   {Cartan derivative}%
\Index{texAffine}
   {Cartan symbol}%
   {Cartan symbol}%
\Index{texAffine}
   {Cartan transport}%
   {Cartan transport}%
\Index{texCartesian}
   {Cartesian power $A$ of set $B$}%
   {Cartesian power of set}%
\Index{texCartesian}
   {Cartesian power of bundle}%
   {Cartesian power of bundle}%
\Index{texCartesian}
   {direct product of bundles}%
   {Cartesian product of bundles}%
\Index{texCartesian}
   {direct product of total spaces}%
   {Cartesian product of total spaces}%
\Index{texHomotopy}
   {category of \drc vector spaces}%
   {category of drc vector spaces}%
\Index{texBundleRelation}
   {category of fibered correspondences over diagonal}%
   {category of fibered correspondences over diagonal}%
\Index{texBundleRelation}
   {category of reduced fibered correspondences}%
   {category of reduced fibered correspondences}%
\Index{texBasis}
   {central affine basis}%
   {Central Affine Basis}%
\Index{texRepresentation}
   {column vector}%
   {column vector}%
\Index{texBundleRelation}
   {commutative diagram of correspondences}%
   {commutative diagram of correspondences}%
\Index{texCartesian}
   {compact\hyph open topology}%
   {compact open topology}%
\Index{texDiffEq}
   {completely integrable system}%
   {completely integrable system}%
\Index{texBundleRelation}
   {composition of fibered correspondences}%
   {composition of fibered correspondences}%
\SubIndex{texVectorSpace}
   {of \drc linear equations}%
   {extended matrix, system of drc linear equations}%
\SubIndex{texVectorSpace}
   {of \rcd linear equations}%
   {extended matrix, system of rcd linear equations}%
\Index{texBundleRelation}
   {composition of reduced fibered correspondences}%
   {composition of reduced fibered correspondences}%
\Index{texBiring}
   {condition of reducibility of products}%
   {condition of reducibility of products}%
\Index{texFiberedGroup}
   {contravariant \Ts representation of fibered group}%
   {Tstar contravariant representation of fibered group}%
\Index{texBasis}
   {coordinate isomorphism}%
   {coordinate isomorphism}%
\Index{texVectorSpace}
   {coordinate isomorphism}%
   {coordinate isomorphism}%
\Index{texVectorSpace}
   {coordinate matrix}%
   {}%
\SubIndex{texVectorSpace}
   {of set of vectors in \dcr rows vector space}%
   {coordinate matrix of set of vectors, dcr vector space}%
\SubIndex{texVectorSpace}
   {of set of vectors in \drc rows vector space}%
   {coordinate matrix of set of vectors, drc vector space}%
\SubIndex{texVectorSpace}
   {of vector in \drc basis}%
   {coordinate matrix of vector in drc basis}%
\Index{texReferenceFrame}
   {coordinate reference frame}%
   {coordinate reference frame}%
\Index{texGenRelativity}
   {coordinate reference frame in event space}%
   {coordinate reference frame in event space}%
\Index{texDrcBasis}
   {coordinate representation in \drc vector space}%
   {coordinate representation, vector space}%
\Index{texBasis}
   {coordinate representation of group in vector space}%
   {coordinate representation, vector space}%
\Index{texBasis}
   {coordinate vector space}%
   {coordinate vector space}%
\Index{texDrcBasis}
   {coordinates of geometrical object}%
   {}%
\SubIndex{texDrcBasis}
   {in coordinate vector space}%
   {coordinates of geometrical object, coordinate vector space}%
\SubIndex{texDrcBasis}
   {in vector space}%
   {coordinates of geometrical object, vector space}%
\Index{texBasis}
   {coordinates of geometrical object}%
   {coordinates of geometrical object, vector space}%
\Index{texBasis}
   {coordinates of geometrical object in coordinate representation}%
   {coordinates of geometrical object, coordinate vector space}%
\Index{texDrcBasis}
   {coordinates of representation}%
   {coordinates of representation}%
\Index{texBasis}
   {coordinates of representation}%
   {coordinates of representation}%
\Index{texVectorSpace}
   {coordinates of set of vectors in \dcr vector space}%
   {coordinates of set of vectors, dcr vector space}%
\Index{texVectorSpace}
   {coordinates of set of vectors in \drc vector space}%
   {coordinates of set of vectors, drc vector space}%
\Index{texVectorSpace}
   {coordinates of vector in \drc basis}%
   {coordinates of vector in drc basis}%
\Index{texBundleRelation}
   {correspondence of homomorphism}%
   {correspondence of homomorphism}%
\Index{texFiberedGroup}
   {covariant \Ts representation of fibered group}%
   {Tstar covariant representation of fibered group}%
\Index{texBiring}
   {\CR inverse element of biring}%
   {cr-inverse element}%
\Index{texVectorSpace}
   {\CR matrix group}%
   {cr-matrix group}%
\Index{texBiring}
   {\CR power}%
   {cr power}%
\Index{texBiring}
   {\CR product of matrices}%
   {cr-product of matrices}%
\Index{texVectorSpace}
   {\crd vector space}%
   {crd vector space}%
\SetIndexSpace%
\Index{texCalculus}
   {$D$\Hyph valued variable}%
   {D valued variable}%
\Index{texVectorSpace}
   {\dcr basis of \subs rows vector space}%
   {dcr basis, c rows vector space}%
\Index{texVectorSpace}
   {\dcr vector}%
   {dcr vector}%
\Index{texVectorSpace}
   {\dcr vector space}%
   {dcr vector space}%
\Index{texBiring}
   {determinant of matrix}%
   {determinant}%
\Index{texTidal}
   {deviation of trajectories}%
   {deviation of trajectories}%
\Index{texBundleRelation}
   {diagonal in bundle}%
   {diagonal in bundle}%
\Index{texBundleRelation}
   {diagram of correspondences}%
   {diagram of correspondences}%
\Index{texCalculus}
   {differentiable functions of \drc vector space to skew field $D$}%
   {differentiable functions, drc vector space to skew field}%
\Index{texCalculus}
   {differential of mapping of normed \drc vector space to valued skew field}%
   {differential, drc vector space to skew field}%
\Index{texVectorSpace}
   {dimension of \drc vector space}%
   {dimension of vector space}%
\Index{texFiberedGroup}
   {direct product of representations of fibered group}%
   {direct product of representations of fibered group}%
\Index{texRepresentation}
   {direct product of representations of group}%
   {direct product of representations of group}%
\Index{texTstarRepresentation}
   {direct product of \Ts representations of group}%
   {direct product of representations of group}%
\Index{texLieRepresentation}
   {direct sum of representations}%
   {direct sum of representations}%
\Index{texVectorSpace}
   {distributive law}%
   {}%
\SubIndex{texVectorSpace}
   {\Ts vector space}%
   {distributive law, Tstar vector space}%
\Index{texVectorSpace}
   {\drc automorphism of vector space}%
   {automorphism of vector space}%
\Index{texVectorSpace}
   {\drc basis}%
   {}%
\SubIndex{texVectorSpace}
   {for \sups rows vector space}%
   {drc basis, r rows vector space}%
\SubIndex{texVectorSpace}
   {for vector space}%
   {drc basis, vector space}%
\Index{texVectorSpace}
   {\drc coordinate vector space}%
   {drc coordinate vector space}%
\Index{texVectorSpace}
   {\drc isomorphism of vector spaces}%
   {isomorphism of vector spaces}%
\Index{texDrcMorphism}
   {\drc linear map of vector spaces}%
   {drc linear map of vector spaces}%
\Index{texVectorSpace}
   {\drc linear span in vector space}%
   {linear span, vector space}%
\Index{texVectorSpace}
   {\drc linearly dependent vectors}%
   {linearly dependent, vector space}%
\Index{texReferenceFrame}
   {\drc linearly independent vector fields}%
   {linearly independent vector fields}%
\Index{texVectorSpace}
   {\drc linearly independent vectors}%
   {linearly independent, vector space}%
\Index{texVectorSpace}
   {\drc system of linear equations}%
   {system of linear equations}%
\Index{texVectorSpace}
   {\drc vector}%
   {drc vector}%
\Index{texCalculus}
   {\drc vector function}%
   {drc vector function}%
\Index{texVectorSpace}
   {\drc vector space}%
   {drc vector space}%
\Index{texVectorSpace}
   {$D\star$\hyph  vector space}%
   {Dstar vector space}%
\Index{texVectorSpace}
   {$D\star$\hyph product of vector over scalar}%
   {Dstar product of vector over scalar, vector space}%
\Index{texBiring}
   {duality principle for biring}%
   {duality principle for biring}%
\Index{texBiring}
   {duality principle for biring of matrices}%
   {duality principle for biring of matrices}%
\SetIndexSpace%
\Index{texTstarMorphism}
   {effective representation of algebra $A$}%
   {effective representation of algebra}%
\Index{texFiberedAlgebra}
   {effective representation of fibered $\mathcal{F}$\Hyph algebra}%
   {effective representation of fibered F-algebra}%
\Index{texFiberedGroup}
   {effective \Ts representation of fibered group}%
   {effective representation of fibered group}%
\Index{texRepresentation}
   {effective representation of group}%
   {effective representation of group}%
\Index{texVectorSpace}
   {effective representation of skew field}%
   {effective representation of skew field}%
\Index{texELie}
   {enhanced Lie group}%
   {enhanced Lie group}%
\Index{texDiffEq}
   {essential parameters}%
   {essential parameters}%
\Index{texBundleRelation}
   {extension of correspondence}%
   {extension of correspondence}%
\Index{texAffine}
   {extreme line}%
   {extreme line}%
\SetIndexSpace%
\Index{texBundleRelation}
   {fibered correspondence from $\mathcal{A}$ to $\mathcal{B}$}%
   {fibered correspondence from A to B}%
\Index{texBundleRelation}
   {fibered correspondence in $\mathcal{A}$}%
   {fibered correspondence in A}%
\Index{texBundleRelation}
   {fibered correspondence of homomorphism}%
   {fibered correspondence of homomorphism}%
\Index{texBundleRelation}
   {fibered equivalence}%
   {fibered equivalence}%
\Index{texFiberedAlgebra}
   {fibered $\mathcal{F}$\Hyph algebra}%
   {fibered F-algebra}%
\Index{texFiberedAlgebra}
   {fibered $\mathcal{F}$\Hyph subalgebra}%
   {fibered F-subalgebra}%
\Index{texFiberedAlgebra}
   {fibered group}%
   {fibered group}%
\Index{texBundleRelation}
   {fibered ordering}%
   {fibered ordering}%
\Index{texBundleRelation}
   {fibered preordering}%
   {fibered preordering}%
\Index{texFiberedAlgebra}
   {fibered ring}%
   {fibered ring}%
\Index{texBundleRelation}
   {fibered subset}%
   {fibered subset}%
\Index{texNewton}
   {field-strength tensor}%
   {field-strength tensor}%
\Index{texNewton}
   {first Newton law}%
   {First Newton law}%
\Index{texAffine}
   {Frenet transport}%
   {Frenet transport}%
\Index{texCalculus}
   {function continuous with respect to set of arguments}%
   {function continuous with respect to set of arguments}%
\Index{texCalculus}
   {function of $\gi n$ $D$\Hyph valued variables}%
   {function of n D valued variables}%
\SetIndexSpace%
\Index{texTypeBasis}
   {\Gbasis}%
   {G-basis}%
\Index{texBasis}
   {\Gbasis\ of vector space}%
   {G-basis}%
\Index{texBasis}
   {\Gcoords\ of basis}%
   {G-coordinates}%
\Index{texTypeBasis}
   {\Gcoords}%
   {G-coordinates}%
\Index{texTypeBasis}
   {\Gspace}%
   {GSpace}%
\Index{texDrcBasis}
   {geometrical object}%
   {}%
\SubIndex{texDrcBasis}
   {defined in vector space}%
   {geometrical object, vector space}%
\SubIndex{texDrcBasis}
   {in coordinate representation defined in vector space}%
   {geometrical object, coordinate vector space}%
\SubIndex{texDrcBasis}
   {of type $A$}%
   {geometrical object of type A, vector space}%
\Index{texBasis}
   {geometrical object in coordinate representation}%
   {geometrical object, coordinate vector space}%
\Index{texBasis}
   {geometrical object in vector space}%
   {geometrical object, vector space}%
\Index{texBasis}
   {geometrical object of type $A$ in vector space}%
   {geometrical object of type A, vector space}%
\Index{texGroupRing}
   {group algebra}%
   {group algebra}%
\Index{texBasis}
   {\Gspace}%
   {GSpace}%
\SetIndexSpace%
\Index{texBiring}
   {Hadamard inverse of matrix}%
   {Hadamard inverse of matrix}%
\Index{texReferenceFrame}
   {holonomic coordinates of vector}%
   {vector holonomic coordinates}%
\Index{texFiberedGroup}
   {homogeneous bundle of fibered group}%
   {homogeneous bundle of fibered group}%
\Index{texTstarRepresentation}
   {homogeneous space of group}%
   {homogeneous space of group}%
\Index{texRepresentation}
   {homogeneous space of group}%
   {homogeneous space of group}%
\Index{texFiberedAlgebra}
   {homomorphism of fibered $\mathcal{F}$\Hyph algebras}%
   {homomorphism of fibered F-algebras}%
\Index{texFiberedGroup}
   {homomorphism of fibered groups}%
   {homomorphism of fibered groups}%
\SetIndexSpace%
\Index{texLieRepresentation}
   {infinitesimal generator}%
   {infinitesimal generator}%
\Index{texLinearLie}
   {infinitesimal generators of group Lie}%
   {infinitesimal generators of group Lie}%
\Index{texDrcBasis}
   {invariance principle}%
   {invariance principle}%
\Index{texBasis}
   {invariance principle in vector space}%
   {invariance principle, vector space}%
\Index{texBundleRelation}
   {inverse fibered correspondence}%
   {inverse fibered correspondence}%
\Index{texBundleRelation}
   {inverse reduced fibered correspondence}%
   {inverse reduced fibered correspondence}%
\Index{texFiberedAlgebra}
   {isomorphism of fibered $\mathcal{F}$\Hyph algebras}%
   {isomorphism of fibered F-algebras}%
\SetIndexSpace%
\Index{texFiberedGroup}
   {kernel of inefficiency of representation of fibered group}%
   {kernel of inefficiency of representation of fibered group}%
\Index{texRepresentation}
   {kernel of inefficiency of representation of group}%
   {kernel of inefficiency of representation of group}%
\Index{texTstarRepresentation}
   {kernel of inefficiency of \Ts representation of group $G$}%
   {kernel of inefficiency of representation of group}%
\Index{texDiffProperty}
   {Killing equation}%
   {Killing equation}%
\Index{texDiffProperty}
   {Killing equation of second type}%
   {Killing equation second type}%
\Index{texDiffProperty}
   {Killing vector of second type}%
   {Killing vector second type}%
\Index{texBiring}
   {Kronecker symbol}%
   {Kronecker symbol}%
\SetIndexSpace%
\Index{texLie}
   {left invariant vector field}%
   {left invariant vector}%
\Index{texVectorSpace}
   {left module}%
   {left module}%
\Index{texTstarRepresentation}
   {left shift}%
   {left shift}%
\Index{texFiberedGroup}
   {left shift on fiberd group}%
   {Tstar shift, fibered group}%
\Index{texRepresentation}
   {left shift on group}%
   {left shift, group}%
\Index{texLie}
   {left structural constant of Lie algebra}%
   {left structural constant of Lie algebra}%
\Index{texVectorSpace}
   {left vector space}%
   {left vector space}%
\Index{texRepresentation}
   {left-side contravariant representation of group}%
   {left-side contravariant representation}%
\Index{texRepresentation}
   {left-side covariant representation of group}%
   {left-side covariant representation}%
\Index{texTstarMorphism}
   {left-side representation of $\mathcal{F}$\Hyph algebra $A$ in set $M$}%
   {left-side representation of algebra}%
\Index{texFiberedAlgebra}
   {left-side representation of fibered $\mathcal{F}$\Hyph algebra}%
   {left-side representation of fibered F-algebra}%
\Index{texRepresentation}
   {left-side representation of group}%
   {left-side representation of group}%
\Index{texTstarMorphism}
   {left-side transformation}%
   {left-side transformation}%
\Index{texFiberedAlgebra}
   {left-side transformation on bundle}%
   {left-side transformation of bundle}%
\Index{texLie}
   {Lie algebra of Lie group}%
   {algebra Lie group Lie}%
\SubIndex{texLie}
   {left defined}%
   {left defined algebra Lie}%
\SubIndex{texLie}
   {right defined}%
   {right defined algebra Lie}%
\Index{texDiffProperty}
   {Lie derivative}%
   {Lie derivative}%
\SubIndex{texDiffProperty}
   {of connection}%
   {Lie derivative of connection}%
\SubIndex{texDiffProperty}
   {of metric}%
   {Lie derivative of metric}%
\Index{texLie}
   {Lie group basic operators}%
   {Lie group basic operators}%
\Index{texBundleRelation}
   {lift of correspondence}%
   {lift of correspondence}%
\Index{texCartesian}
   {lift of map}%
   {lift of map}%
\Index{texRepresentation}
   {linear representation of group}%
   {linear representation of group}%
\Index{texReferenceFrame}
   {linearly dependent vector fields}%
   {linearly dependent vector fields}%
\Index{texReferenceFrame}
   {local reference frame}%
   {local reference frame}%
\Index{texGenRelativity}
   {Lorentz transformation}%
   {Lorentz transformation}%
\SetIndexSpace%
\Index{texReferenceFrame}
   {map of type $G$ on manifold}%
   {map of type G on manifold}%
\Index{texCartesian}
   {mapping space}%
   {mapping space}%
\Index{texDrcMorphism}
   {matrix of \drc linear map}%
   {matrix of drc linear map}%
\Index{texAffine}
   {metric-affine manifold}%
   {metric-affine manifold}%
\Index{texTstarMorphism}
   {morphism of \Ts representations from $f$ into $g$}%
   {morphism of representations from f into g}%
\Index{texTstarMorphism}
   {morphism of \Ts representations of $\mathcal{F}$ algebra}%
   {morphism of representations of F algebra}%
\Index{texTstarMorphism}
   {morphism of \Ts representations of $\mathcal{F}$\Hyph algebra in $\mathcal{H}$\Hyph algebra}%
   {morphism of representations of F algebra in H algebra}%
\Index{}
   {movement on basis manifold}%
   {movement transformation}%
\SetIndexSpace%
\Index{texPolymodule}
   {$(n)$\hyph vector space}%
   {(n)-vector space}%
\Index{texBundleRelation}
   {$n$\Hyph ary fibered relation}%
   {fibered relation}%
\Index{texAffine}
   {nonmetricity}%
   {nonmetricity}%
\Index{texVectorSpace}
   {nonsingular \drc system of linear equations}%
   {nonsingular system of linear equations}%
\Index{texRepresentation}
   {nonsingular \Ts transformation}%
   {nonsingular transformation}%
\Index{texCalculus}
   {norm on \drc vector space}%
   {norm on drc vector space}%
\Index{texCalculus}
   {normed \drc vector space}%
   {normed drc vector space}%
\SetIndexSpace%
\Index{texFiberedAlgebra}
   {operation on bundle}%
   {operation on bundle}%
\Index{texBundleRelation}
   {opposite fibered preordering}%
   {opposite fibered preordering}%
\Index{texFiberedGroup}
   {orbit of representation of fibered group}%
   {orbit of representation of fibered group}%
\Index{texRepresentation}
   {orbit of representation of group}%
   {orbit of representation of group}%
\Index{texTstarRepresentation}
   {orbit of \Ts representation of group}%
   {orbit of representation of group}%
\Index{texBasis}
   {orthonornal basis}%
   {Orthonornal Basis}%
\SetIndexSpace%
\Index{texReferenceFrame}
   {parallelogram}%
   {parallelogram}%
\Index{texCalculus}
   {partial derivative of mapping $f$ with respect to variable $v^{\gi a}$}%
   {partial derivative of mapping with respect to variable, skew field}%
\Index{texCalculus}
   {partial derivative of mapping $\Vector f$ with respect to variable $v^{\gi a}$}%
   {partial derivative of mapping with respect to variable, drc vector space}%
\Index{texBasis}
   {passive representation}%
   {passive representation}%
\Index{texBasis}
   {passive transformation on basis manifold}%
   {passive transformation}%
\Index{texDrcBasis}
   {passive \Ts representation}%
   {passive representation}%
\Index{texReferenceFrame}
   {pfaffian derivative}%
   {pfaffian derivative}%
\Index{texNewton}
   {potential energy}%
   {potential energy}%
\Index{texDrcBasis}
   {product of geometrical object and constant}%
   {product of geometrical object and constant}%
\Index{texBasis}
   {product of geometrical object and constant in vector space}%
   {product of geometrical object and constant, vector space}%
\Index{texTstarMorphism}
   {product of morphisms of \Ts representations of $\mathcal{F}$\Hyph algebra}%
   {product of morphisms of representations of F algebra}%
\SetIndexSpace%
\Index{texBasis}
   {quasi affine transformation on basis manifold}%
   {quasi affine transformation}%
\Index{texBasis}
   {quasi movement on basis manifold}%
   {quasi movement}%
\SetIndexSpace%
\Index{texBiring}
   {\sups row of matrix}%
   {r row}%
\Index{texVectorSpace}
   {\sups rows \drc vector space}%
   {sups rows drc vector space}%
\Index{texBiring}
   {$r$\hyph row of matrix}%
   {r-row}%
\Index{texBiring}
   {\RC inverse element of biring}%
   {rc-inverse element}%
\Index{texVectorSpace}
   {\RC major minor}%
   {RC-major minor}%
\Index{texVectorSpace}
   {\RC matrix group}%
   {rc-matrix group}%
\Index{texVectorSpace}
   {\RC nonsingular matrix}%
   {RC nonsingular matrix}%
\Index{texBiring}
   {\RC power}%
   {rc power}%
\Index{texBiring}
   {\RC product of matrices}%
   {rc-product of matrices}%
\Index{texBiring}
   {\RC quasideterminant}%
   {RC-quasideterminant}%
\Index{texVectorSpace}
   {\RC rank of matrix}%
   {rc-rank of matrix}%
\Index{texVectorSpace}
   {\RC singular matrix}%
   {RC singular matrix}%
\Index{texVectorSpace}
   {\rcd vector space}%
   {rcd vector space}%
\Index{texCartesian}
   {reduced Cartesian product of bundles}%
   {reduced Cartesian product of bundles}%
\Index{texCartesian}
   {reduced Cartesian product of total spaces}%
   {reduced Cartesian product of total spaces}%
\Index{texBundleRelation,texBundleRelation}
   {reduced fibered correspondence from $\mathcal{A}$ to $\mathcal{B}$}%
   {reduced fibered correspondence from A to B}%
\Index{texBundleRelation}
   {reduced fibered correspondence in $\mathcal{A}$}%
   {reduced fibered correspondence in A}%
\Index{texBiring}
   {reducible biring}%
   {reducible biring}%
\Index{texReferenceFrame}
   {reference frame}%
   {reference frame}%
\Index{texGenRelativity}
   {reference frame in event space}%
   {reference frame in event space}%
\Index{texReferenceFrame}
   {reference frame manifold}%
   {reference frame manifold}%
\Index{texBundleRelation}
   {reflexive $2$\Hyph ary fibered relation}%
   {reflexive 2 ary fibered relation}%
\Index{texTstarRepresentation}
   {representation of group}%
   {}%
\SubIndex{texTstarRepresentation}
   {contravariant \Ts}%
   {Tstar contravariant representation of group}%
\SubIndex{texTstarRepresentation}
   {covariant \Ts}%
   {Tstar covariant representation of group}%
\SubIndex{texDrcBasis}
   {\drc linear \sT}%
   {linear representation of group}%
\SubIndex{texTstarRepresentation}
   {effective}%
   {effective representation of group}%
\SubIndex{texDrcBasis}
   {\rcd}%
   {rcd linear representation of group}%
\SubIndex{texTstarRepresentation}
   {\sT}%
   {starT representation of group}%
\SubIndex{texTstarRepresentation}
   {\Ts}%
   {Tstar representation of group}%
\Index{texRepresentation}
   {representation of group}%
   {representation of group}%
\Index{texDrcBasis}
   {representative of geometrical object in vector space}%
   {representative of geometrical object, vector space}%
\Index{texBasis}
   {representative of geometrical object in vector space}%
   {representative of geometrical object, vector space}%
\Index{texBundleRelation}
   {restriction of correspondence $\Phi$ to set $C$}%
   {restriction of correspondence}%
\Index{texLie}
   {right invariant vector field}%
   {right invariant vector}%
\Index{texVectorSpace}
   {right module}%
   {right module}%
\Index{texTstarRepresentation}
   {right shift}%
   {right shift}%
\Index{texRepresentation}
   {right shift on group}%
   {right shift, group}%
\Index{texLie}
   {right structural constant of Lie algebra}%
   {right structural constant of Lie algebra}%
\Index{texVectorSpace}
   {right vector space}%
   {right vector space}%
\Index{texRepresentation}
   {right-side contravariant representation of group}%
   {right-side contravariant representation}%
\Index{texRepresentation}
   {right-side covariant representation of group}%
   {right-side covariant representation}%
\Index{texTstarMorphism}
   {right-side representation of $\mathcal{F}$\Hyph algebra $A$ in set $M$}%
   {right-side representation of algebra}%
\Index{texFiberedAlgebra}
   {right-side representation of fibered $\mathcal{F}$\Hyph algebra}%
   {right-side representation of fibered F-algebra}%
\Index{texRepresentation}
   {right-side representation of group}%
   {right-side representation of group}%
\Index{texTstarMorphism}
   {right-side transformation}%
   {right-side transformation}%
\Index{texRepresentation}
   {right-side transformation}%
   {right-side transformation}%
\Index{texRepresentation}
   {row vector}%
   {row vector}%
\Index{texVectorSpace}
   {$R\star$\Hyph module}%
   {Rstar-module}%
\SetIndexSpace%
\Index{texNewton}
   {scalar potential}%
   {scalar potential}%
\Index{texNewton}
   {second Newton law}%
   {Second Newton law}%
\Index{texTstarMorphism}
   {single transitive representation of algebra $A$}%
   {single transitive representation of algebra}%
\Index{texFiberedAlgebra}
   {single transitive representation of fibered $\mathcal{F}$\Hyph algebra}%
   {single transitive representation of fibered F-algebra}%
\Index{texRepresentation}
   {single transitive representation of group}%
   {single transitive representation of group}%
\Index{texTidal}
   {speed of deviation}%
   {speed of deviation}%
\Index{texDrcMorphism}
   {$(S\RCstar,T\RCstar)$\Hyph linear map of vector spaces}%
   {src trc linear map of vector spaces}%
\Index{texDrcBasis}
   {standard coordinates of basis}%
   {standard coordinates of basis}%
\Index{texBasis}
   {standard coordinates of basis}%
   {standard coordinates of basis}%
\Index{texBiring}
   {standard representation of matrix}%
   {Standard representation}%
\Index{texVectorSpace}
   {$\star D$\hyph vector space}%
   {starD-vector space}%
\Index{texLinearMap}
   {$\star D$\Hyph product of \drc linear map $A$ over scalar}%
   {starD product of drc linear map over scalar}%
\Index{texVectorSpace}
   {$\star R$\hyph module}%
   {starR-module}%
\Index{texTstarRepresentation}
   {\sT shift}%
   {starT shift}%
\Index{texFiberedGroup}
   {\sT shift on fibered group}%
   {starT shift, fibered group}%
\Index{texTstarMorphism}
   {\sT representation of $\mathcal{F}$\Hyph algebra $A$ in set $M$}%
   {starT representation of algebra}%
\Index{texFiberedAlgebra}
   {\sT representation of fibered $\mathcal{F}$\Hyph algebra}%
   {starT representation of fibered F-algebra}%
\Index{texFiberedGroup}
   {\sT representation of fibered group}%
   {starT representation of fibered group}%
\Index{texTstarMorphism}
   {\sT transformation}%
   {starT transformation}%
\Index{texFiberedAlgebra}
   {\sT transformation on bundle}%
   {starT transformation of bundle}%
\SubIndex{}
   {nonsingular}%
   {nonsingular transformation of bundle}%
\Index{texBundleRelation}
   {subbundle}%
   {subbundle}%
\Index{texLinearMap}
   {sum of \drc linear maps}%
   {sum of drc linear maps, drc vector spaces}%
\Index{texDrcBasis}
   {sum of geometrical objects}%
   {sum of geometrical objects}%
\Index{texBasis}
   {sum of geometrical objects in vector space}%
   {sum of geometrical objects, vector space}%
\Index{texBundleRelation}
   {symmetric $2$\Hyph ary fibered relation}%
   {symmetric 2 ary fibered relation}%
\Index{texBasis}
   {symmetry group}%
   {symmetry group}%
\Index{texDrcBasis}
   {symmetry group}%
   {SymmetryGroup}%
\Index{texGenRelativity}
   {synchronization of reference frame}%
   {synchronization of reference frame}%
\SetIndexSpace%
\Index{texLie}
   {tensor product of representations}%
   {tensor product of representations}%
\Index{texCalculus}
   {topological \drc vector space}%
   {topological drc vector space}%
\Index{texCalculus}
   {topological skew field}%
   {topological skew field}%
\Index{texAffine}
   {torsion}%
   {torsion}%
\Index{texTstarMorphism}
   {transformation of set}%
   {transformation of set}%
\Index{texDrcBasis}
   {transformation on basis manifold}%
   {}%
\SubIndex{texDrcBasis}
   {active}%
   {active transformation, vector space}%
\SubIndex{texTypeBasis}
   {affine}%
   {affine transformation}%
\SubIndex{texTypeBasis}
   {movement}%
   {movement transformation}%
\SubIndex{texDrcBasis}
   {passive}%
   {passive transformation, vector space}%
\SubIndex{texTypeBasis}
   {quasi affine}%
   {quasi affine transformation}%
\SubIndex{texTypeBasis}
   {quasi movement}%
   {quasi movement}%
\Index{texFiberedAlgebra}
   {transformation on bundle}%
   {transformation of bundle}%
\Index{texBundleRelation}
   {transitive $2$\Hyph ary fibered relation}%
   {transitive 2 ary fibered relation}%
\Index{texTstarMorphism}
   {transitive representation of algebra $A$}%
   {transitive representation of algebra}%
\Index{texFiberedAlgebra}
   {transitive representation of fibered $\mathcal{F}$\Hyph algebra}%
   {transitive representation of fibered F-algebra}%
\Index{texRepresentation}
   {transitive representation of group}%
   {transitive representation of group}%
\Index{texVectorSpace}
   {\Ts linear composition of  vectors}%
   {linear composition of  vectors}%
\Index{texVectorSpace}
   {\Ts matrices vector space}%
   {matrices vector space}%
\Index{texTstarRepresentation}
   {\Ts shift}%
   {Tstar shift}%
\Index{texTstarMorphism}
   {\Ts representation of $\mathcal{F}$\Hyph algebra $A$ in set $M$}%
   {Tstar representation of algebra}%
\Index{texFiberedAlgebra}
   {\Ts representation of fibered $\mathcal{F}$\Hyph algebra}%
   {Tstar representation of fibered F-algebra}%
\Index{texFiberedGroup}
   {\Ts representation of fibered group}%
   {Tstar representation of fibered group}%
\Index{texTstarMorphism}
   {\Ts transformation}%
   {Tstar transformation}%
\Index{texFiberedAlgebra}
   {\Ts transformation on bundle}%
   {Tstar transformation of bundle}%
\Index{texFiberedGroup}
   {twin representations of fibered group}%
   {twin representations of fibered group}%
\Index{texTstarRepresentation}
   {twin representations of group}%
   {twin representations of group}%
\Index{texLinearMap}
   {twin representations of skew field}%
   {twin representations of skew field}%
\Index{texReferenceFrame}
   {type $G$ reference frame}%
   {type G reference frame}%
\SetIndexSpace%
\Index{texVectorSpace}
   {unitarity law}%
   {}%
\SubIndex{texVectorSpace}
   {for \Ts vector space}%
   {unitarity law, Tstar vector space}%
\SetIndexSpace%
\Index{texPolymodule}
   {($D_1\RCstar$, ..., $D_n\RCstar$)\hyph vector space}%
   {(d1rc,dnrc)-vector space}%
\Index{texPolymodule}
   {($S\star$, $\star T$)\hyph vector space}%
   {(Sstar,starT)-vector space}%
\Index{texCalculus}
   {valued skew field}%
   {valued skew field}%
\Index{texFiberedAlgebra}
   {vector bundle}%
   {vector bundle}%
\Index{texNewton}
   {vector potential}%
   {vector potential}%
\Index{texVectorSpace}
   {vector space type}%
   {vector space type}%

\CloseIndex

\def\indexname{Special Symbols and Notations}
\OpenIndex

\SetIndexSpace%
\Symb{texBiring}
   {$(^{\gi a}_{\gi b})$\hyph\CR quasideterminant}%
   {a b CR quasideterminant definition}%
\Symb{texBiring}
   {$(^{\gi a}_{\gi b})$\hyph \RC quasideterminant}%
   {a b RC-quasideterminant definition}%
\Symb{texBiring}
   {minor}%
   {A from b a}%
\Symb{texBiring}
   {minor}%
   {A from columns T}%
\Symb{texBiring}
   {minor}%
   {A from rows S}%
\Symb{texBiring}
   {minor}%
   {A without column a}%
\Symb{texBiring}
   {minor}%
   {A without columns T}%
\Symb{texBiring}
   {minor}%
   {A without row b}%
\Symb{texBiring}
   {minor}%
   {A without rows S}%
\Symb{texPolymodule}
   {active transformation}%
   {active transformation}%
\Symb{texTypeBasis}
   {affine space}%
   {affine space}%
\Symb{texBasis}
   {affine space}%
   {An}%
\Symb{texBiring}
   {\subs row ($c$\hyph row) of matrix}%
   {c row}%
\Symb{texBiring}
   {\CR power of element $A$ of biring}%
   {cr power}%
\Symb{texBiring}
   {\CR inverse element of biring}%
   {cr-inverse element}%
\Symb{texBiring}
   {\CR product of matrices}%
   {cr-product of matrices}%
\Symb{texVectorSpace}
   {\dcr vector}%
   {dcr vector}%
\Symb{texLie}
   {derivative of left shift}%
   {derivative of left shift}%
\Symb{texLie}
   {derivative of left shift}%
   {derivative of left shift, 1-Parameter Group}%
\Symb{texLie}
   {derivative of right shift}%
   {derivative of right shift}%
\Symb{texLie}
   {}%
   {derivative of right shift}%
\Symb{texLie}
   {derivative of right shift}%
   {derivative of right shift, 1-Parameter Group}%
\Symb{texLie}
   {derivative of left shift}%
   {derivative of Tstar shift}%
\Symb{texVectorSpace}
   {\drc vector}%
   {drc vector}%
\Symb{texAffine}
   {derivative}%
   {overline nabla_l, definition 2}%
\Symb{texPolymodule}
   {passive transformation}%
   {passive transformation}%
\Symb{texBiring}
   {\sups row ($r$\hyph row) of matrix}%
   {r row}%
\Symb{texBiring}
   {\RC power of element $A$ of biring}%
   {rc power}%
\Symb{texBiring}
   {\RC inverse element of biring}%
   {rc-inverse element}%
\Symb{texBiring}
   {\RC product of matrices}%
   {rc-product of matrices}%
\Symb{texBiring}
   {\RC quasideterminant}%
   {RC-quasideterminant definition}%
\Symb{texTstarRepresentation}
   {right shift}%
   {right shift}%
\Symb{texFiberedGroup}
   {\sT shift}%
   {starT shift, fibered group}%
\Symb{texTstarRepresentation}
   {left shift}%
   {Tstar shift}%
\Symb{texFiberedGroup}
   {\Ts shift}%
   {Tstar shift, fibered group}%
\Symb{texReferenceFrame}
   {anholonomic coordinates of vector}%
   {vector anholonomic coordinates}%
\Symb{texReferenceFrame}
   {holonomic coordinates of vector}%
   {vector holonomic coordinates}%

\SetIndexSpace%
\Symb{texBasis}
   {basis manifold of affine space}%
   {BAn}%
\Symb{texReferenceFrame}
   {basis manifold of manifold}%
   {basis manifold of manifold}%
\Symb{texPolymodule}
   {basis manifold of vector space}%
   {basis manifold of vector space}%
\Symb{texBasis}
   {basis manifold of vector space $\mathcal{V}$}%
   {basis manifold of vector space}%
\Symb{texBasis}
   {basis manifold of central affine space}%
   {BCAn}%
\Symb{texBasis}
   {basis manifold of Euclid space}%
   {BEn}%
\Symb{texCartesian}
   {Cartesian power $A$ of set $B$}%
   {Cartesian power of set}%
\Symb{texTypeBasis}
   {basis manifold of affine space}%
   {FAn}%
\Symb{texTypeBasis}
   {basis manifold of central affine space}%
   {FCAn}%
\Symb{texTypeBasis}
   {basis manifold of Euclid space}%
   {FEn}%

\SetIndexSpace%
\Symb{texBasis}
   {central affine space}%
   {CAn}%
\Symb{texTypeBasis}
   {central affine space}%
   {central affine space}%
\Symb{texLie}
   {left structural constant of Lie algebra}%
   {left structural constant of Lie algebra}%
\Symb{texLie}
   {right structural constant of Lie algebra}%
   {right structural constant of Lie algebra}%
\Symb{texReferenceFrame}
   {set of functions defined on manifold}%
   {set of functions defined on manifold}%

\SetIndexSpace%
\Symb{texLieRepresentation}
   {basis vector of \sT representation}%
   {basis vector of starT representation}%
\Symb{texLieRepresentation}
   {basis vector of \sT representation}%
   {basis vector of starT representation, coordinates}%
\Symb{texLieRepresentation}
   {basis vector of \Ts representation}%
   {basis vector of Tstar representation}%
\Symb{texLieRepresentation}
   {basis vector of \Ts representation}%
   {basis vector of Tstar representation, coordinates}%
\Symb{texVectorSpace}
   {\subs rows \dcr vector space}%
   {c rows dcr vector space}%
\Symb{texCalculus}
   {differential of function}%
   {differential, drc vector space to drc vector space}%
\Symb{texCalculus}
   {differential of function}%
   {differential, drc vector space to skew field}%
\Symb{texVectorSpace}
   {\drc coordinate vector space}%
   {drc coordinate vector space}%
\Symb{texVectorSpace}
   {matrices vector space}%
   {matrices vector space}%
\Symb{texAffine}
   {Cartan derivative}%
   {overbrace D}%
\Symb{texAffine}
   {derivative}%
   {overline D}%
\Symb{texCalculus}
   {partial derivative of mapping $\Vector f$ with respect to variable $v^{\gi a}$}%
   {partial derivative of mapping, 1, drc vector space}%
\Symb{texCalculus}
   {partial derivative of mapping $f$ with respect to variable $v^{\gi a}$}%
   {partial derivative of mapping, 1, skew field}%
\Symb{texVectorSpace}
   {\sups rows \drc vector space}%
   {r rows drc vector space}%
\Symb{texTidal}
   {speed of deviation}%
   {speed of deviation}%
\Symb{texVectorSpace}
   {vector space type}%
   {vector space type}%

\SetIndexSpace%
\Symb{texTypeBasis}
   {affine basis}%
   {Affine Basis}%
\Symb{texBasis}
   {affine basis}%
   {Affine Basis}%
\Symb{texTypeBasis}
   {basis}%
   {basis}%
\Symb{texBasis}
   {basis of vector space}%
   {Basis e}%
\Symb{texBasis}
   {basis in vector space $\mathcal{V}$}%
   {basis in V}%
\Symb{texVectorSpace}
   {basis of vector space}%
   {basis, vector space}%
\Symb{texPolymodule}
   {basis of $(n)$\hyph vector space}%
   {basis,n vector space}%
\Symb{texCartesian}
   {Cartesian power of total spaces}%
   {Cartesian power of total spaces}%
\Symb{texCartesian}
   {Cartesian product of total spaces}%
   {Cartesian product of total spaces, definition 1}%
\Symb{texBasis}
   {central affine basis}%
   {Central Affine Basis}%
\Symb{texReferenceFrame}
   {form of reference frame}%
   {dual forms, reference frame}%
\Symb{texTypeBasis}
   {Euclid space}%
   {En}%
\Symb{texBasis}
   {Euclid space}%
   {En}%
\Symb{texBasis}
   {pseudo Euclid space}%
   {Enm}%
\Symb{texTypeBasis}
   {pseudo Euclid space}%
   {Enm}%
\Symb{texFiberedAlgebra}
   {identical transformation of bundle}%
   {identical transformation of bundle}%
\Symb{texBasis}
   {orthonornal basis}%
   {Orthonornal Basis}%
\Symb{texCartesian}
   {reduced Cartesian product of total spaces}%
   {reduced Cartesian product of total spaces, definition 1}%
\Symb{texFiberedAlgebra}
   {set of nonsingular \sT transformations of bundle $\mathcal{E}$}%
   {set of starT nonsingular transformations of bundle}%
\Symb{texFiberedAlgebra}
   {set of nonsingular \Ts transformations of bundle $\mathcal{E}$}%
   {set of Tstar nonsingular transformations of bundle}%
\Symb{texBasis}
   {standard coordinates of basis}%
   {standard coordinates of basis}%
\Symb{texReferenceFrame}
   {standard coordinates of reference frame}%
   {standard coordinates of reference frame}%
\Symb{texReferenceFrame}
   {vector field of reference frame}%
   {vector field of reference frame}%
\Symb{texReferenceFrame}
   {vector field of reference frame}%
   {vector field, reference frame}%
\Symb{texBasis}
   {vector of basis}%
   {vector of basis}%

\SetIndexSpace%
\Symb{texVectorSpace}
   {coordinates of basis in \subs rows \dcr vector space}%
   {basis coordinates, c rows dcr vector space}%
\Symb{texVectorSpace}
   {coordinates of basis in \sups rows \drc vector space}%
   {basis coordinates, r rows drc vector space}%
\Symb{texVectorSpace}
   {basis for \subs rows \dcr vector space}%
   {basis, c rows dcr vector space}%
\Symb{texVectorSpace}
   {basis for \sups rows \drc vector space}%
   {basis, r rows drc vector space}%
\Symb{texDiffEq}
   {central affine basis}%
   {Central Affine Basis}%
\Symb{texReferenceFrame}
   {coordinate reference frame}%
   {coordinate reference frame}%
\Symb{texFiberedAlgebra}
   {homomorphism of fibered $\mathcal{F}$\Hyph algebras}%
   {homomorphism of fibered F-algebras}%
\Symb{texBundleRelation}
   {inverse fibered correspondence}%
   {inverse fibered correspondence, 1}%
\Symb{texBundleRelation}
   {inverse reduced fibered correspondence}%
   {inverse reduced fibered correspondence, 1}%
\Symb{texCartesian}
   {map to Cartesian product}%
   {map to Cartesian product}%
\Symb{texTstarRepresentation}
   {representation orbit of group $G$}%
   {orbit of representation of group}%
\Symb{texTypeBasis}
   {orthonornal basis}%
   {Orthonornal Basis}%
\Symb{texReferenceFrame}
   {reference frame}%
   {reference frame}%
\Symb{texReferenceFrame}
   {reference frame}%
   {reference frame, extensive definition}%
\Symb{texPolymodule}
   {standard coordinates of basis}%
   {standard coordinates of basis}%
\Symb{texPolymodule}
   {vector of basis}%
   {vector of basis}%

\SetIndexSpace%
\Symb{texVectorSpace}
   {\CR matrix group}%
   {cr-matrix group}%
\Symb{texLie}
   {algebra Lie of group Lie}%
   {g}%
\Symb{texLie}
   {left defined algebra Lie of group Lie}%
   {gl}%
\Symb{texTypeBasis}
   {affine transformation group}%
   {GLAn}%
\Symb{texBasis}
   {affine transformation group}%
   {GLAn}%
\Symb{texLie}
   {right defined algebra Lie of group Lie}%
   {gr}%
\Symb{texBasis}
   {group of homomorphisms of vector space $\mathcal{V}$}%
   {GV}%
\Symb{texFiberedGroup}
   {orbit of effective covariant \sT representation of fibered group}%
   {orbit of effective starT covariant representation of fibered group}%
\Symb{texTstarRepresentation}
   {orbit of effective covariant \sT representation of group}%
   {orbit of effective starT covariant representation of group}%
\Symb{texFiberedGroup}
   {orbit of effective covariant \Ts representation of fibered group}%
   {orbit of effective Tstar covariant representation of fibered group}%
\Symb{texTstarRepresentation}
   {orbit of effective covariant \Ts representation of group}%
   {orbit of effective Tstar covariant representation of group}%
\Symb{texVectorSpace}
   {\RC matrix group}%
   {rc-matrix group}%

\SetIndexSpace%
\Symb{texBiring}
   {Hadamard inverse of matrix}%
   {Hadamard inverse of matrix}%
\Symb{texLinearMap}
   {\rcd vector space of \drc linear maps}%
   {rcd vector space of drc linear maps}%

\SetIndexSpace%
\Symb{texLieRepresentation}
   {infinitesimal generator of representation}%
   {infinitesimal generator of representation}%
\Symb{texLinearLie}
   {Lie group infinitesimal generators}%
   {Lie group infinitesimal generators}%

\SetIndexSpace%
\Symb{texRepresentation}
   {left shift}%
   {left shift}%
\Symb{texDiffProperty}
   {Lie derivative of connection}%
   {Lie derivative of connection}%
\Symb{texDiffProperty}
   {Lie derivative of metric}%
   {Lie derivative of metric}%
\Symb{texBasis}
   {passive transformation}%
   {passive transformation}%
\Symb{texRepresentation}
   {set of left-side nonsingular transformations of set $M$}%
   {set of left-side nonsingular transformations}%

\SetIndexSpace%
\Symb{texTstarMorphism}
   {set of \sT transformations of set $M$}%
   {set of starT transformations}%
\Symb{texTstarMorphism}
   {set of \Ts transformations of set $M$}%
   {set of Tstar transformations}%

\SetIndexSpace%
\Symb{texBasis}
   {geometrical object in coordinate representation}%
   {geometrical object, coordinate vector space}%
\Symb{texBasis}
   {geometrical object}%
   {geometrical object, vector space}%
\Symb{texFiberedGroup}
   {orbit of representation of fibered group $\mathcal{G}$}%
   {orbit of representation of fibered group}%
\Symb{texRepresentation}
   {orbit of representation of the group $G$}%
   {orbit of representation of group}%

\SetIndexSpace%
\Symb{texCartesian}
   {bundle}%
   {bundle}%
\Symb{texCartesian}
   {Cartesian power of bundle}%
   {Cartesian power of bundle}%
\Symb{texCartesian}
   {Cartesian product of bundles}%
   {Cartesian product of bundles, definition 1}%
\Symb{texCartesian}
   {reduced Cartesian product of bundles}%
   {reduced Cartesian product of bundles, definition 1}%
\Symb{texFiberedAlgebra}
   {set of nonsingular \sT transformations of bundle $\bundle{}pE{}$}%
   {set of starT nonsingular transformations of bundle, projection}%
\Symb{texFiberedAlgebra}
   {set of nonsingular \Ts transformations of bundle $\bundle{}pE{}$}%
   {set of Tstar nonsingular transformations of bundle, projection}%

\SetIndexSpace%
\Symb{texBasis}
   {active transformation}%
   {active transformation}%
\Symb{texAffine}
   {Cartan curvature}%
   {Cartan curvature}%
\Symb{texVectorSpace}
   {\CR rank of matrix}%
   {cr-rank of matrix}%
\Symb{texBundleRelation}
   {diagonal in bundle  $\bundle{}pA{}$}%
   {diagonal in bundle, 2}%
\Symb{texBundleRelation}
   {diagonal in bundle $\mathcal{A}$}%
   {diagonal in reduced bundle, 2}%
\Symb{texAffine}
   {curvature}%
   {GLn curvature_overline}%
\Symb{texVectorSpace}
   {\RC rank of matrix}%
   {rc-rank of matrix}%
\Symb{texRepresentation}
   {right shift}%
   {right shift}%
\Symb{texRepresentation}
   {set of right-side nonsingular transformations of set $M$}%
   {set of right-side nonsingular transformations}%

\SetIndexSpace%
\Symb{texBundleRelation}
   {composition of fibered correspondences}%
   {composition of fibered correspondences}%
\Symb{texBundleRelation}
   {inverse fibered correspondence}%
   {inverse fibered correspondence, 2}%
\Symb{texBundleRelation}
   {inverse reduced fibered correspondence}%
   {inverse reduced fibered correspondence, 2}%
\Symb{texVectorSpace}
   {linear span in vector space}%
   {linear span, vector space}%

\SetIndexSpace%
\Symb{texLie}
   {tangent plane to group $G$}%
   {TaG}%

\SetIndexSpace%
\Symb{texBasis}
   {coordinate vector space}%
   {coordinate vector space}%
\Symb{texBasis}
   {coordinates in vector space}%
   {coordinates in vector space}%
\Symb{texVectorSpace}
   {\dcr vector space}%
   {left CR vector space}%
\Symb{texVectorSpace}
   {\drc vector space}%
   {left RC vector space}%
\Symb{texLinearMap}
   {($S$, $T$)\hyph bimodule}%
   {R S bimodule}%
\Symb{texVectorSpace}
   {\crd vector space}%
   {right CR vector space}%
\Symb{texVectorSpace}
   {\rcd vector space}%
   {right RC vector space}%
\Symb{texBasis}
   {vector space}%
   {V}%
\Symb{texReferenceFrame}
   {vector space of vector fields}%
   {vector space of vector fields}%

\SetIndexSpace%
\Symb{texPolymodule}
   {geometrical object in coordinate representation		defined in vector space}%
   {geometrical object, coordinate vector space}%
\Symb{texPolymodule}
   {geometrical object in vector space}%
   {geometrical object, vector space}%

\SetIndexSpace%
\Symb{texReferenceFrame}
   {anholonomic coordinate}%
   {x(k)}%

\SetIndexSpace%
\Symb{texBundleRelation}
   {diagonal in bundle $\mathcal{A}$}%
   {diagonal in bundle, 1}%

\SetIndexSpace%
\Symb{texTidal}
   {deviation of trajectories}%
   {deviation of trajectories}%
\Symb{texRepresentation}
   {identical transformation}%
   {identical transformation}%
\Symb{texTstarMorphism}
   {identical transformation}%
   {identical transformation}%
\Symb{texBasis}
   {image of vector $\Vector e_k\in\Basis e$ under isomorphism to coordinate vector space}%
   {image of vector e_k, coordinate vector space}%
\Symb{texBiring}
   {Kronecker symbol}%
   {Kronecker symbol}%

\SetIndexSpace%
\Symb{texReferenceFrame}
   {anholonomic coordinates of connection}%
   {anholonomic coordinates of connection}%
\Symb{texAffine}
   {Cartan symbol}%
   {Cartan symbol}%
\Symb{texAffine}
   {connection}%
   {conection overline}%
\Symb{texAffine}
   {Cartan connection}%
   {overbrace Gamma i kl}%
\Symb{texCartesian}
   {set of sections of bundle}%
   {set of sections of bundle}%

\SetIndexSpace%
\Symb{texLie}
   {inverse operator to operator $\psi_l$}%
   {inverse operator to operator psi l}%
\Symb{texLie}
   {inverse operator to operator $\psi_r$}%
   {inverse operator to operator psi r}%

\SetIndexSpace%
\Symb{texReferenceFrame}
   {anholonomity object}%
   {anholonomity object}%

\SetIndexSpace%
\Symb{texLie}
   {basic operator of group Lie}%
   {Lie Basic Operator L}%
\Symb{texLie}
   {}%
   {Lie Basic Operator L}%
\Symb{texLie}
   {basic operator of group Lie}%
   {Lie Basic Operator L, 1-Parameter Group}%
\Symb{texLie}
   {basic operator of group Lie}%
   {Lie Basic Operator R}%
\Symb{texLie}
   {}%
   {Lie Basic Operator R}%
\Symb{texLie}
   {basic operator of group Lie}%
   {Lie Basic Operator R, 1-Parameter Group}%

\SetIndexSpace%
\Symb{texReferenceFrame}
   {coordinate reference frame}%
   {coordinate reference frame, extensive definition}%
\Symb{texCalculus}
   {partial derivative of mapping $\Vector f$ with respect to variable $v^{\gi a}$}%
   {partial derivative of mapping, 2, drc vector space}%
\Symb{texCalculus}
   {partial derivative of mapping $f$ with respect to variable $v^{\gi a}$}%
   {partial derivative of mapping, 2, skew field}%
\Symb{texReferenceFrame}
   {derivative $e_{(k)}$}%
   {partial(k)}%

\SetIndexSpace%
\Symb{texLie}
   {Lie group composition law}%
   {Lie group composition law}%

\SetIndexSpace%
\Symb{texAffine}
   {}%
   {overbrace nabla_l}%
\Symb{texAffine}
   {Cartan derivative}%
   {overbrace nabla_l}%
\Symb{texAffine}
   {derivative}%
   {overline nabla_l, definition 1}%

\SetIndexSpace%
\Symb{texBundleRelation}
   {restriction of correspondence $\Phi$ to set $C$}%
   {restriction of correspondence}%

\SetIndexSpace%
\Symb{texCartesian}
   {Cartesian product of bundles}%
   {Cartesian product of bundles, definition 2}%
\Symb{texCartesian}
   {Cartesian product of total spaces}%
   {Cartesian product of total spaces, definition 2}%
\Symb{texCartesian}
   {reduced Cartesian product of bundles}%
   {reduced Cartesian product of bundles, definition 2}%
\Symb{texCartesian}
   {reduced Cartesian product of total spaces}%
   {reduced Cartesian product of total spaces, definition 2}%

\SetIndexSpace%
\Symb{texBundleRelation}
   {fibered subset}%
   {fibered subset}%
\Symb{texBundleRelation}
   {subbundle}%
   {subbundle}%

\CloseIndex

\end{document}